\documentclass[a4paper, 10pt]{article}

\usepackage[top=1.5in,bottom=1.5in,right=1.5in,left=1.5in]{geometry}
\usepackage{amsfonts}
\usepackage{amsmath}
\usepackage{amssymb}
\usepackage{amsthm}
\usepackage{fancyref}
\usepackage{mathrsfs}
\usepackage{enumitem}
\usepackage{color}
\usepackage{soul}
\usepackage[normalem]{ulem}
\usepackage{romannum}
\usepackage{mathtools}
\usepackage{csquotes}
\usepackage[UKenglish]{babel}
\usepackage[UKenglish]{isodate}
\usepackage[capitalise]{cleveref}
  \crefformat{equation}{\textup{#2(#1)#3}}
\usepackage{cite}

\theoremstyle{plain}
\newtheorem{theorem}{Theorem}[section]
\newtheorem{lemma}[theorem]{Lemma}
\newtheorem{proposition}[theorem]{Proposition}

\theoremstyle{remark}
\newtheorem{remark}[theorem]{Remark}
\theoremstyle{definition}
\newtheorem{definition}[theorem]{Definition}

\newtheorem{assumption}[theorem]{Assumption}

\crefname{assumption}{Assumption}{Assumptions}
\Crefname{assumption}{Assumption}{Assumptions}

\newcounter{counter_a}
\newenvironment{myenum}{\begin{list}{\textrm{\textup{(\roman{counter_a})}}}%
{\usecounter{counter_a}
\setlength{\itemsep}{0.5ex}\setlength{\topsep}{0.7ex}
\setlength{\leftmargin}{5ex}\setlength{\labelwidth}{5ex}}}{\end{list}}

\newcounter{counter_b}
\newenvironment{myenuma}{\begin{list}{\textrm{\textup{(\alph{counter_b})}}}%
{\usecounter{counter_b}
\setlength{\itemsep}{0.5ex}\setlength{\topsep}{0.7ex}
\setlength{\leftmargin}{5ex}\setlength{\labelwidth}{5ex}}}{\end{list}}

\newenvironment{mybull}{\begin{list}{$\bullet$}%
{\setlength{\itemsep}{0.5ex}\setlength{\topsep}{0.7ex}
\setlength{\leftmargin}{5ex}\setlength{\labelwidth}{5ex}}}{\end{list}}

\newcommand\rd{\mathrm{d}}
\newcommand\NN{\mathbb{N}}
\newcommand\RR{\mathbb{R}}
\newcommand\cB{\mathcal{B}}
\newcommand\cD{\mathcal{D}}
\newcommand\cI{\mathcal{I}}

\newcommand\cL{\mathcal{L}}
\newcommand\cM{\mathcal{M}}
\newcommand\mr{\mathring}
\newcommand\normcdot{\lVert\,\cdot\,\rVert}
\newcommand\normcdotsub[1]{\lVert\,\cdot\,\rVert_{#1}}

\newcommand{\restrict}{\raisebox{-0.25ex}{$\big|$}}
\newcommand{\hatM}{\rule{0ex}{1ex}\hspace{0.45ex}\widehat{\rule[1.5ex]{1.5ex}{0ex}}\hspace{-2.2ex}M}
\newcommand{\tildew}{\widetilde{w}}
\newcommand{\tmax}{t_{\textup{\textsf{max}}}}
\newcommand{\totmass}{\phi_{[1]}}
\newcommand{\totmassM}{\cM_1}
\definecolor{grey}{rgb}{0.4,0.4,0.4}

\numberwithin{equation}{section}

\begin{document}

\pagenumbering{arabic}

\title{Discrete Coagulation--Fragmentation Systems in \\[0.5ex] Weighted $\ell^1$ Spaces}
\author{Lyndsay Kerr\footnote{Department of Mathematics and Statistics,
University of Strathclyde,
26 Richmond Street,
Glasgow G1 1XH,
United Kingdom;
email: \texttt{lyndsay.kerr@strath.ac.uk}, \texttt{m.langer@strath.ac.uk}; \\
ORCID: 0000-0002-6667-7175, 0000-0001-8813-7914} 
\; and Matthias Langer$^*$}
\date{}

\maketitle

\begin{abstract}
\noindent
We study an infinite system of ordinary differential equations that models the evolution of 
coagulating and fragmenting clusters, which we assume to be composed of identical units.
Under very mild assumptions on the coefficients we prove existence, uniqueness and positivity 
of solutions of a corresponding semi-linear Cauchy problem in a weighted $\ell^1$ space.
This requires the application of novel results, which we prove for abstract semi-linear Cauchy problems 
in Banach lattices where the non-linear term is defined only on a dense subspace.
\\[1ex]
\textit{Mathematics Subject Classification \textup{(}2020\textup{)}:}  34G20; 47D06, 80A30, 46B42, 46B70
\\[1ex]
\textit{Keywords:} discrete coagulation--fragmentation, analytic semigroup, semi-linear abstract Cauchy problem
\end{abstract}


\section{Introduction}\label{C01}

We consider a system consisting of clusters of particles that can merge together
(i.e.\ coagulate) to produce larger clusters and can break apart (i.e.\ fragment)
to produce smaller clusters.
Assume that each cluster is made up of identical units, monomers, and that
the mass is scaled such that a monomer has unit mass.
Then a cluster comprised of $n \in \mathbb{N}$ monomers (an $n$-mer) has mass $n$.
We also say that an $n$-mer has size $n$.
Hence cluster size and mass are discrete variables, the values of which coincide.
If we take the number density of clusters of size $n$ at time $t\in[0,T)$ to be $u_n(t)$,
where $0<T\le\infty$, then the evolution of $n$-mers can be described by the
coagulation--fragmentation (C--F) equations
\begin{equation}\label{C21}
\begin{split}
  u_n'(t) &= {}-a_nu_n(t)+\sum\limits_{j=n+1}^{\infty} a_jb_{n,j}u_j(t) \\
  &\quad +\frac{1}{2}\sum\limits_{j=1}^{n-1} k_{n-j,j}(t)u_{n-j}(t)u_j(t)
  -\sum\limits_{j=1}^{\infty} k_{n,j}(t)u_n(t)u_j(t), \qquad t \in (0,T); \\[1ex]
  u_n(0) &= \mathring{u}_n,
\end{split}
\end{equation}
$n=1,2,\ldots$.
The first two terms in \eqref{C21} represent the change in the density
of clusters of size $n$ due to fragmentation, while the final two terms describe
the change in the density of clusters of size $n$ due to coagulation.
For $n \ge 2$, $a_n$ is the rate at which clusters of size $n \in \mathbb{N}$ fragment.
Since monomers cannot fragment, when $a_1>0$, we interpret $a_1$ to be the rate at
which monomers are removed from the system through some other mechanism.
The fragmentation coefficient $b_{n,j}$, for $n<j$, is the average number
of $n$-mers that are produced when a $j$-mer fragments.
By $k_{n,j}(t)$, we denote the rate at which $n$-mers merge with $j$-mers at time $t$.
Note that we allow the coagulation coefficients, $k_{n,j}$, to be time-dependent, 
whereas the fragmentation coefficients, $a_n$, $b_{n,j}$, are independent of time.
Finally, $\mathring{u}_n$ denotes the initial density of clusters of size $n$ at time $t=0$.
We can express the solution of \eqref{C21},
as a time-dependent sequence $u(t)=(u_n(t))_{n=1}^{\infty}$.
Since, for all $t \ge 0$, $u(t)$ is a sequence of densities, we expect $u(t)$
to be a non-negative sequence.

We consider the following non-negativity assumption on the rate coefficients,
which we assume to hold throughout this paper.

\begin{assumption}\label{C22}
Let $0<T\le\infty$.
Assume that, for all $n,j\in \mathbb{N}$, $a_n \ge 0$, $b_{n,j} \ge 0$ and $k_{n,j}(t) \ge 0$
for all $t\in[0,T)$.  Moreover, $b_{n,j}=0$ for all $j \le n$ and $k_{n,j}(t)=k_{j,n}(t)$
for all $n, j \in \mathbb{N}$, $t\in[0,T)$.
\hfill $\lozenge$
\end{assumption}

\medskip

\noindent
If $u(t)$ is a non-negative solution of the system \eqref{C21},
then the total mass of clusters at time $t \ge 0$ is given by
\begin{equation}\label{C23}
	\totmassM(t) \coloneq \sum_{n=1}^{\infty} nu_n(t).
\end{equation}
When a cluster of size $j$ fragments, the total mass of the
daughter clusters that are produced is equal to
\begin{equation}\label{C24}
	\sum\limits_{n=1}^{j-1} nb_{n,j} 
\end{equation}
for $j\ge2$.
A formal calculation shows that one expects
\begin{equation}\label{C25}
	\frac{\rd}{\rd t}\totmassM(t) = \sum\limits_{n=1}^{\infty} nu_n'(t)
	= -a_1u_1(t) - \sum_{j=2}^\infty \biggl(j-\sum_{n=1}^{j-1}nb_{n,j}\biggr)a_ju_j(t)
\end{equation}
for all $t\in(0,T)$.
Hence, if
\begin{equation}\label{C149}
	\sum_{n=1}^{j-1}nb_{n,j} \le j \qquad\text{for all} \ j\in\{2,3,\ldots\},
\end{equation}
then the mass is expected to be non-increasing;
if 
\begin{equation}\label{C26}
	a_1=0 \qquad \text{and} \qquad
	\sum\limits_{n=1}^{j-1} nb_{n,j}=j \qquad \text{for all} \ j\in\{2,3,\ldots\},
\end{equation}
then the mass is expected to be conserved.
However, in contrast to most papers in the literature, we do not assume in our main results
that \eqref{C149} or \eqref{C26} hold, i.e.\ we allow that mass is gained or lost
during fragmentation events.

Processes of coagulation and fragmentation occur in a diverse range of areas
including polymer science,
\cite{aizenman_bak1979, ziff1980kinetics, ziff_mcgrady1985kinetics},
in the formation of aerosols, \cite{drake1972aerosol},
and in the powder production industry, \cite{verdurmen_etal2004simulation, wells2018thesis}.
The system \eqref{C21} originates from work by Smoluchowski,
who derived an equation to describe pure coagulation in 1916/1917,
\cite{smoluchowski1916coagulation,smoluchowski1917coagulation}.
In 1935 Becker and D\"oring studied a special case with both coagulation and fragmentation \cite{becker_doring1935}.

In general, neither existence nor uniqueness of solutions of \eqref{C21} are guaranteed;
for non-uniqueness for the pure fragmentation equation, see, for instance, the discussion
in \cite[Example~4.3]{kerr_lamb_langer2020fragpaper}.
One approach that is used to deal with \eqref{C21} is a weak compactness argument.
In this approach the full C--F system is truncated and a sequence of solutions to
these truncated equations is obtained.
A weak compactness argument can be used to show that this sequence converges
to a solution of an integrated version of the C--F system,
and further work then enables the uniqueness of this solution to be shown.
This approach was first used by Ball, Carr and Penrose, in \cite{ball_carr_penrose1986becker},
to deal with the Becker--D\"oring cluster equations and has since been used to
deal with the case of binary fragmentation (where particles can fragment only
into two smaller particles) in, e.g.\ \cite{ball_carr1990discrete, dacosta1995existence},
and to deal with the case of multiple fragmentation (where particles can fragment
into several smaller particles) in, e.g.\ \cite{laurencot2002discretemultiple}; see also \cite{banasiak_lamb_laurencot2020vol2}.

In this paper we use a semigroup approach to deal with \eqref{C21}. 
This approach was pioneered in \cite{aizenman_bak1979} to examine the binary fragmentation version of the continuous C--F system, 
and has since been used extensively to deal with both the discrete and continuous system 
\cite{banasiak2001extension, mcbride_smith_lamb2010strongly, mclaughlin_lamb_mcbride1997fragmentation, mclaughlin_lamb_mcbride1997CF, smith_lamb_langer_mcbride2012discrete};
see also \cite[Section~2.2.1]{banasiak_lamb_laurencot2020vol1} and the references therein.
We begin by formulating \eqref{C21} as a
semi-linear ACP (abstract Cauchy problem) in an appropriate Banach space.
The fragmentation terms in \eqref{C21} give rise to a linear operator
in this Banach space, which is known to be the generator of a substochastic $C_0$-semigroup;
see, e.g.\ \cite{kerr_lamb_langer2020fragpaper}.
We then treat the coagulation operator, originating from the coagulation terms
in \eqref{C21}, as a non-linear perturbation and apply
perturbation theory to obtain results regarding the existence
and uniqueness of solutions.

Since we can write the solution of \eqref{C21}
as $u(t)=\bigl(u_n(t)\bigr)_{n=1}^{\infty}$, we work in an infinite sequence space.
It turns out that weighted $\ell^1$ spaces are most appropriate.
Let the sequence $w=(w_n)_{n=1}^{\infty}$ be such that $w_n>0$ for all $n\in\NN$.
We reformulate \eqref{C21} as a semi-linear ACP in the space
\begin{equation}\label{C96}
	\ell_w^1 \coloneq \biggl\{f=(f_n)_{n=1}^\infty:
	f_n \in \mathbb{R} \ \text{for all} \ n \in \mathbb{N} \ \text{and} \
	\sum\limits_{n=1}^{\infty} w_n|f_n|<\infty\biggr\},
\end{equation}
equipped with the norm
\[
	\|f\|_w \coloneq \sum_{n=1}^{\infty} w_n|f_n| \qquad \text{for all} \ f \in \ell_w^1.
\]
We refer to this Banach space as the weighted $\ell^1$ space with weight $w$.

Motivated by \eqref{C23}, it is most common to work with $w_n=n$
for all $n\in\NN$; see \cite{mcbride_smith_lamb2010strongly}.
In this particular case we set $X_{[1]} \coloneqq \ell_w^1$, i.e.\
\begin{equation}\label{C147}
	X_{[1]} \coloneq \biggl\{f=(f_n)_{n=1}^\infty:
	f_n\in\RR \ \text{for all} \ n \in \mathbb{N} \ \text{and} \
	\sum\limits_{n=1}^{\infty} n|f_n|<\infty\biggr\},
\end{equation}
and $\|f\|_{[1]}=\sum_{n=1}^\infty n|f_n|$.
This is the most physically relevant Banach space to work in as the norm,
$\normcdotsub{[1]}$, of a non-negative solution in this space coincides with
the total mass in the system. 
Let us introduce the linear functional
\begin{equation}\label{C159}
	\totmass(f) \coloneq \sum_{n=1}^\infty nf_n, \qquad f=(f_j)_{j=1}^\infty \in X_{[1]},
\end{equation}
which coincides with $\normcdotsub{[1]}$ on the set of positive elements.
Then the total mass in \eqref{C23} of a solution $u(t)$ can be written in terms of $\totmass$ 
as $\totmassM(t)=\totmass(u(t))$.
Some investigation has also been undertaken into the case where $w_n=n^p$ for $p>1$;
see, e.g.\ \cite{banasiak2012global}.
When $w_n=n^p$, $n\in\NN$, $p\ge1$, we set $X_{[p]} \coloneqq \ell_w^1$
and $\normcdotsub{[p]} \coloneqq \normcdotsub{w}$.
In the present paper we work with more general weights and follow on from \cite{kerr_lamb_langer2020fragpaper} 
where the linear part of \eqref{C21}, i.e.\ the pure fragmentation equation without coagulation, was studied.
As demonstrated in that paper \cite{kerr_lamb_langer2020fragpaper}, working with a more general weight
can be very beneficial as it allows us to utilise results about analytic semigroups
for a larger class of fragmentation coefficients.
In particular, analytic semigroups enable us to use interpolation spaces.
The latter have been used for continuous and discrete coagulation--fragmentation equations in the literature; 
see, e.g.\ \cite{banasiak2012global,banasiak_lamb2012analytic,banasiak_lamb_langer2013strong}.
However, with the flexibility of choosing the weight we can use interpolation spaces even in cases of very slow fragmentation.

In \cref{C100}, the main result of our paper, we prove local existence, uniqueness and positivity of mild solutions of 
the coagulation--fragmentation equation \eqref{C21} under very mild assumptions on the coagulation coefficients, $k_{n,j}$.
In particular, the flexibility of choosing the weight $w$ enables us to allow very general coagulation coefficients.
Let us also point out that, in prior investigations of \eqref{C21}, it has been assumed that
the coagulation rates are independent of time, i.e.\ $k_{n,j}(t)\equiv k_{n,j}$
for all $n,j\in\NN$, $t \ge 0$;
see, e.g.\ \cite{banasiak2012global} and \cite{mcbride_smith_lamb2010strongly}.
Here we relax this assumption and allow time-dependent coagulation.

In order to prove \cref{C100}, we establish results about existence and uniqueness of solutions
of abstract semi-linear ACPs in a Banach space (or Banach lattice) $X$; see \cref{C31,C32}.
The novel features of these theorems is the inclusion of positivity of solutions under weak assumptions,
which is important both for physical applications and to prove other results,
and the use of a space $Y$, continuously embedded in $X$, such that the non-linear term
is defined only on $Y$.
The flexibility of the space $Y$ allows us to apply the abstract theorems in different situations,
including the case when $Y$ is an interpolation space.

Let us give a brief summary of the contents of the paper.
We begin in Section~\ref{C02} by providing some preliminary concepts and results, 
in particular relating to Lipschitz continuity and Fr\'echet differentiability. 
We also supply useful results relating to operator semigroups and interpolation spaces. 
In Section~\ref{C03} we then provide the abstract results relating to the existence and uniqueness 
of solutions to semi-linear ACPs, formulated in terms of general Banach spaces $X$ and $Y$, 
with $Y$ continuously embedded in $X$.
The proofs of \cref{C31,C32} are given in \cref{C16,C17} respectively, whereas
in \cref{C15} we provide sufficient conditions for the assumptions of \cref{C31,C32}
and results about the behaviour of certain linear functionals applied to solutions,
which can be used to show mass conservation or estimates for the total mass.
Finally, in \cref{C04} we apply the abstract results from \cref{C03} in two different settings, 
corresponding to different spaces $Y$, to prove our main result, \cref{C100}, 
about existence, uniqueness and positivity of mild solutions of a semi-linear ACP describing the C--F system.
We finish the paper with results about classical solutions in \cref{C158} and
about global existence of the solution in \cref{C152}.

\section{Preliminaries}\label{C02}

In this section we introduce some terminology and results that we apply later to
the coagulation--fragmentation system.
First, let us fix some notation.
Let $X$ be a normed vector space with norm $\normcdot$.
We define the \emph{open ball} (respectively \emph{closed ball}) centred at $h \in X$
and with radius $r>0$ by
\begin{equation*}
  B_X(h,r) \coloneqq \{g \in X: \|h-g\|<r\} \quad (\text{respectively} \
  \overline{B}_X(h,r) \coloneqq \{g \in X: \Vert h-g \Vert \le r\}).
\end{equation*}
Further, for normed vector spaces $X$ and $Y$, we denote by $\cB(X,Y)$ the space of
all bounded linear operators from $X$ to $Y$,
define $\cB(X)\coloneq\cB(X,X)$, and denote
by $\cD(A)$ the domain of a linear operator $A$.

\subsection{Lipschitz continuity and Fr\'echet differentiability}
\label{C11}

To analyse the coagulation--fragmentation system, we need a local Lipschitz condition
that is locally uniform in the time variable.
The following definition is taken from \cite[p.~185]{pazy1983semigroups}.

\begin{definition}\label{C27}
Let $t_0,T \in \mathbb{R}$ be such that $t_0<T$ and let $\cI$ be an interval of the
form $[t_0,T)$, $[t_0,T]$ or $[t_0,\infty)$.
Further, let $(Y,\normcdotsub{Y})$, $(X,\normcdotsub{X})$ be normed vector spaces.
Then an operator $F:\cI \times Y \to X$ satisfies
a \emph{Lipschitz condition in the second argument on bounded sets, uniformly
in the first argument on compact intervals} if, for every $t' \in \cI$ and constant $r>0$,
there exists a constant $L=L(t',r)$ such that
\begin{equation}\label{C28}
  \|F(t,f)-F(t,g)\|_X \le L \|f-g\|_Y
\end{equation}
for all $f,g \in Y$ with $\|f\|_Y,\|g\|_Y \le r$ and $t \in [t_0,t']$.
\end{definition}

Note that when $\cI=[t_0,T]$, we can take $t'=T$ so that \eqref{C28}
holds with the same $L$ for all $t \in [t_0,T]$.

The following lemma is useful in proving the validity of a Lipschitz condition in the situation of a mapping 
that is quadratic in the space variable;
cf., e.g.\ \cite[Theorems~4.3 and 4.4]{mcbride_smith_lamb2010strongly} for a special case in connection with a coagulation term.

\begin{lemma}\label{C29}
Let $(Y,\normcdotsub{Y})$, $(X,\normcdotsub{X})$ be Banach spaces.
Further, let $t_0$, $T \in \mathbb{R}$ be such that $t_0<T<\infty$ and
let $\cI$ be an interval of the form $[t_0,T)$, $[t_0,T]$ or $[t_0,\infty)$.
Let $\widetilde{F}: \cI \times Y \times Y \to X$ be such that
\begin{myenuma}
\item
  $\widetilde{F}$ is linear in the second and third arguments;
\item
  for each $t' \in \cI$ there exists a $c(t')>0$ such that
  \begin{equation}\label{C30}
    \big\|\widetilde{F}[t,f,g]\big\|_X \le c(t')\|f\|_Y\|g\|_Y
  \end{equation}
  for all $f,g \in Y$ and $t \in [t_0,t']$.
\end{myenuma}
Define $F: \cI \times Y \to X$ by $F(t,f) \coloneqq \widetilde{F}[t,f,f]$
for $t\in\cI$, $f\in Y$.
\begin{myenum}
\item
  Then $F$ is Lipschitz in the second argument on bounded sets, uniformly
  in the first argument on compact intervals.
\item
  If, in addition, $t \mapsto \widetilde{F}[t,f,g]$ is continuous on $\cI$ for
  fixed $f$, $g \in Y$, then the mapping $\widetilde{F}: \cI \times Y \times Y \to X$
  is continuous and hence $F: \cI \times Y \to X$ is continuous.
\end{myenum}
\end{lemma}

\begin{proof}
(i)
Let $t'\in\cI$, $r>0$ and let $f,g \in Y$ be such that $\|f\|_Y,\|g\|_Y \le r$.
Then, for all $t \in [t_0,t']$, we have
\begin{align*}
  \big\|F(t,f)-F(t,g)\big\|_X
  &= \big\|\widetilde{F}[t,f,f]-\widetilde{F}[t,g,f]
  +\widetilde{F}[t,g,f]-\widetilde{F}[t,g,g]\big\|_X \\[1ex]
  &\le \big\|\widetilde{F}[t,f-g,f]\big\|_X + \big\|\widetilde{F}[t,g,f-g]\big\|_X \\[1ex]
  &\le c(t')\|f-g\|_Y \|f\|_Y + c(t')\|g\|_Y \|f-g\|_Y \\[1ex]
  &\le 2rc(t')\|f-g\|_Y.
\end{align*}

(ii)
Now suppose that, in addition, $t \mapsto \widetilde{F}[t,f,g]$ is continuous on $\cI$
for every $f,g\in Y$.  Fix $(s_0,f_0,g_0) \in \cI \times Y \times Y$.
Choose $t'\in\cI$ such that $t'>s_0$ (if $\cI=[t_0,T]$ and $s_0=T$, we take $t'=s_0=T$).
Then, for $s\in[t_0,t']$, $f,g\in Y$,
\begin{align*}
  & \big\|\widetilde{F}(s,f,g)-\widetilde{F}(s_0,f_0,g_0)\big\|_X
  \\[1ex]
  &= \big\|\widetilde{F}(s,f,g)-\widetilde{F}(s,f_0,g)
  +\widetilde{F}(s,f_0,g)-\widetilde{F}(s,f_0,g_0)
  +\widetilde{F}(s,f_0,g_0)-\widetilde{F}(s_0,f_0,g_0)\big\|_X
  \\[1ex]
  &\le \big\|\widetilde{F}(s,f-f_0,g)\big\|_X + \big\|\widetilde{F}(s,f_0,g-g_0)\big\|_X
  + \big\|\widetilde{F}(s,f_0,g_0)-\widetilde{F}(s_0,f_0,g_0)\big\|_X
  \displaybreak[0]\\[1ex]
  &\le c(t')\|f-f_0\|_Y\|g\|_Y + c(t')\|f_0\|_Y\|g-g_0\|_Y
  + \big\|\widetilde{F}(s,f_0,g_0)-\widetilde{F}(s_0,f_0,g_0)\big\|_X
  \\[1ex]
  &\to 0 \qquad \text{as} \;\; (s,f,g) \to (s_0,f_0,g_0).
\end{align*}
It follows that $\widetilde{F}: \cI \times Y \times Y$ is continuous,
and hence so is $F: \cI \times Y$.
\end{proof}

\begin{remark}\label{C146}
Let $F:\cI\times Y\to X$ satisfy a Lipschitz condition in the second argument 
on bounded sets, uniformly in the first argument on compact intervals.
Further, let $V\in\cB(Y,X)$.
Then $F'(t,f)\coloneq F(t,f)+Vf$ also satisfies a Lipschitz condition in the second argument 
on bounded sets, uniformly in the first argument on compact intervals.
\rule{0ex}{1ex}\hfill $\lozenge$
\end{remark}

\medskip

Let again $(X,\normcdotsub{X})$, $(Y,\normcdotsub{Y})$ be Banach spaces and
define the norm $\normcdotsub{\RR\times Y}$ on the space $\RR\times Y$ by
\begin{equation}\label{C40}
  \|(t,f)\|_{\RR\times Y} \coloneqq |t|+\|f\|_Y.
\end{equation}
Let $t_0,T\in\RR$ be such that $t_0<T$ and let $\cI$ be an interval of the
form $[t_0,T)$, $[t_0,T]$ or $[t_0,\infty)$.  Further, let $F: \cI \times Y \to X$
and take $(t,f) \in \mr{\cI} \times Y$, where $\mr{\cI}$ denotes the interior of $\cI$.
Recall that the operator $F$ is \emph{Fr\'echet differentiable} at $(t,f)$
if there exists a linear operator $DF(t,f) \in \cB(\RR\times Y,X)$ and
a mapping $R((t,f), \cdot): \RR \times Y \to X$ such that,
for all $(\delta_t,\delta_f) \in (\RR \times Y)\setminus \{(0,0)\}$
satisfying $t+\delta_t \in \cI\setminus\{t_0\}$,
\begin{equation*}
  F(t+\delta_t,f+\delta_f)
  = F(t,f)+\bigl[DF(t,f)\bigr](\delta_t,\delta_f) + R\bigl((t,f),(\delta_t,\delta_f)\bigr),
\end{equation*}
where
\[
  \lim_{\|(\delta_t,\delta_f)\|_{\RR \times Y} \to 0}
  \frac{\big\|R\bigl((t,f),(\delta_t,\delta_f)\bigr)\big\|_X}{\|(\delta_t,\delta_f)\|_{\RR\times Y}}
  = 0.
\]
The operator $DF(t,f)$ is the \emph{Fr\'echet derivative} of $F$ at $(t,f)$.
The mapping $F$ is said to be \emph{Fr\'echet differentiable} on an
open subset $D_0 \subset \cI \times Y$ if it is Fr\'echet differentiable
at every $(t,f) \in D_0$.

In the next lemma we provide sufficient conditions for a mapping
to be Fr\'echet differentiable.

\begin{lemma}\label{C41}
Let $(Y,\normcdotsub{Y})$, $(X,\normcdotsub{X})$ be Banach spaces
and let $t_0,T \in \mathbb{R}$ be such that $t_0<T$.
Further, let $\cI$ be an interval of the form $[t_0,T)$, $[t_0,T]$ or $[t_0,\infty)$
and let $F: \cI \times Y \to X$ be a mapping that satisfies
\begin{myenuma}
\item
  $f\mapsto F(t,f)$ is Fr\'echet differentiable for every $t\in\cI$,
  with a uniform remainder term on bounded intervals of $t$ in the sense that,
  for all $t\in\cI$ and $f,\delta_f \in Y$ with $\delta_f \ne 0$, we have
  \begin{equation}\label{C42}
    F(t,f+\delta_f) = F(t,f)+\bigl[D_YF(t,f)\bigr]\delta_f+R_Y\bigl((t,f),\delta_f\bigr),
  \end{equation}
  where for fixed $f$,
  \begin{equation}\label{C43}
    \frac{\|R_Y\bigl((t,f),\delta_f\bigr)\|_X}{\|\delta_f\|_Y} \to 0
    \qquad \text{as} \;\; \|\delta_f\|_Y \to 0,
  \end{equation}
  uniformly in $t$ on compact subintervals of $\cI$,
  and $D_Y F(t,f)\in \cB(Y,X)$ for each $(t,f)\in\cI\times Y$;
\item
  $t \mapsto F(t,f)$ is strongly differentiable for fixed $f \in Y$;
\item
  $(t,f) \mapsto D_YF(t,f)$ and $(t,f) \mapsto \frac{\partial}{\partial t} F(t,f)$
  are continuous;
\item
  $t\mapsto \bigl[D_Y F(t,f)\bigr]\delta_f$ is differentiable for all $f,\delta_f\in Y$
  and, for every $t'\in\mr{\cI}$ there exists $\hat c(f,t')>0$ such that
  \[
    \bigg\|\frac{\partial}{\partial t}\bigl[D_Y F(t,f)\bigr]\delta_f\bigg\|_X \le \hat c(f,t')\|\delta_f\|_Y
    \qquad \text{for} \ t\in(t_0,t'], \; \delta_f\in Y.
  \]
\end{myenuma}
Then $F$ is Fr\'echet differentiable on $\mr{\cI}\times Y$ and
the Fr\'echet derivative, $DF(t,f)$, of $F$ at $(t,f)\in\mr{\cI}\times Y$
is given by
\begin{equation}\label{C44}
  \bigl[DF(t,f)\bigr](s,g) = \frac{\partial}{\partial t}F(t,f)s + D_YF(t,f)g
\end{equation}
for all $s\in\RR$, $t\in Y$.
Moreover, $DF(t,f)$ is continuous with respect to $(t,f)$.
\end{lemma}

\begin{proof}
Assumption~(b) implies that, for all $t\in\mr{\cI}$,
$\delta_t \in \RR\setminus\{0\}$ such that $t+\delta_t \in\mr{\cI}$,
and $f \in Y$, we have that
\begin{equation}\label{C45}
  F(t+\delta_t,f) = F(t,f)+\frac{\partial}{\partial t}F(t,f)\delta_t+R_t\bigl((t,f),\delta_t\bigr)
\end{equation}
where
\begin{equation}\label{C94}
  \frac{\big\|R_t\bigl((t,f),\delta_t\bigr)\big\|_X}{|\delta_t|} \to 0
  \qquad \text{as} \;\; |\delta_t| \to 0.
\end{equation}
Let $t\in\mr{\cI}$, $\delta_t \in \RR\setminus\{0\}$
such that $t+\delta_t\in\mr{\cI}$, $f\in Y$ and $\delta_f\in Y\setminus\{0\}$.
Then from \eqref{C42}, \eqref{C45} and
assumption (d) we obtain that
\begin{align*}
  & F(t+\delta_t,f+\delta_f)-F(t,f)
  \\[1ex]
  &= F(t+\delta_t, f+\delta_f) - F(t+\delta_t,f) + F(t+\delta_t,f) - F(t,f)
  \displaybreak[0]\\[1ex]
  &= \bigl[D_YF(t+\delta_t,f)\bigr]\delta_f + R_Y\bigl((t+\delta_t,f),\delta_f\bigr)
  + \frac{\partial}{\partial t} F(t,f)\delta_t+R_t\bigl((t,f),\delta_t\bigr)
  \displaybreak[0]\\[1ex]
  &= \bigl[D_YF(t,f)\bigr]\delta_f+\int\limits_t^{t+\delta_t}
  \frac{\partial}{\partial \tau} \Bigl(\bigl[D_YF(\tau,f)\bigr]\delta_f\Bigr)\,\rd\tau
  + R_Y\bigl((t+\delta_t,f),\delta_f\bigr)
  \\
  &\quad + \frac{\partial}{\partial t}F(t,f)\delta_t + R_t\bigl((t,f),\delta_t\bigr).
\end{align*}
Now fix $(t,f) \in \mr{\cI}\times Y$.
Then $\bigl[D_YF(t,f)\bigr]\delta_f+\frac{\partial}{\partial t} F(t,f)\delta_t$ is
linear in $(\delta_t, \delta_f)$ and bounded since $D_YF(t,f)$ is linear and bounded.
With the norm $\normcdotsub{\RR\times Y}$ as in \eqref{C40} we have
\begin{align*}
  & \frac{1}{\|(\delta_t,\delta_f)\|_{\RR\times Y}}
  \left\|\int\limits_t^{t+\delta_t}
  \frac{\partial}{\partial\tau}\Bigl(\bigl[D_YF(\tau,f)\bigr]\delta_f\Bigr)\,\rd\tau\right\|_X
  \\[1ex]
  &\le \frac{1}{\|\delta_f\|_Y}\int\limits_t^{t+\delta_t}\,
  \biggl\|\frac{\partial}{\partial t}\Bigl(\bigl[D_YF(\tau,f)\bigr]\delta_f\Bigr)\biggr\|\,\rd\tau
  \to 0
\end{align*}
as $\|(\delta_t,\delta_f)\|_{\RR\times Y} \to 0$ by assumption (d).
From \eqref{C43} we obtain
\[
  \frac{\big\|R_Y\bigl((t+\delta_t,f),\delta_f\bigr)\big\|_X}{\|(\delta_t,\delta_f)\|_{\RR\times Y}}
  \le \frac{\big\|R_Y\bigl((t+\delta_t,f),\delta_f\bigr)\|_X}{\|\delta_f\|_Y}
  \to 0 \qquad \text{as} \;\; \|(\delta_t,\delta_f)\|_{\RR\times Y} \to 0
\]
and similarly from \eqref{C94},
\[
  \frac{\big\|R_t\bigl((t,f),\delta_t\bigr)\|_X}{\|(\delta_t,\delta_f)\|_{(\RR\times Y)}}
  \le \frac{\big\|R_t\bigl((t,f),\delta_t\bigr)\|_X}{|\delta_t|} \to 0
  \qquad \text{as} \;\; \|(\delta_t,\delta_f)\|_{\RR\times Y} \to 0.
\]
It follows that $F$ is Fr\'echet differentiable on $\mr{\cI}\times Y$,
and the Fr\'echet derivative of $F$ at $(t,f)\in\mr{\cI}\times Y$
is given by \eqref{C44}.

Since $D_YF(t,f)$ and $\frac{\partial}{\partial t} F(t,f)$ are continuous
with respect to $(t,f)$, it follows from \eqref{C44}
that $DF(t,f)$ is also continuous with respect to $(t,f)$.
\end{proof}

In the next result we come back to the more particular situation in Lemma~\ref{C29}
and provide sufficient conditions for the Fr\'echet differentiability of such mappings.

\begin{proposition}\label{C46}
Let the conditions of Lemma~\ref{C29}\,\textup{(ii)}
be satisfied with $\cI=[t_0,T)$, $\cI=[t_0,T]$ or $\cI=[t_0,\infty)$.  Moreover, assume that
\begin{myenuma}
\item
  for fixed $f,g \in Y$, $t \mapsto \widetilde{F}[t,f,g]$ is continuously differentiable on $\mr{\cI}$;
\item
  for each $t'\in\cI$, there exists $\tilde{c}(t')>0$ such that
  \begin{equation*}
    \bigg\|\frac{\partial}{\partial t}\widetilde{F}[t,f,g]\bigg\|_X
    \le \tilde{c}(t')\|f\|_Y \|g\|_Y
  \end{equation*}
  for all $f,g\in Y$ and $t\in(t_0,t']$.
\end{myenuma}
Then $F$ is Fr\'echet differentiable with respect to $(t,f)$
on\, $\mr{\cI}\times Y$, and the Fr\'echet derivative
at $(t,f)\in\mr{\cI}\times Y$ is given by
\[
  \bigl[DF(t,f)\bigr](s,g) = \Bigl(\frac{\partial}{\partial t}\widetilde{F}[t,f,f]\Bigr)s
  + \widetilde{F}[t,f,g] + \widetilde{F}[t,g,f],
\]
for all $(s,g)\in\RR\times Y$.  Moreover, $DF(t,f)$ is continuous with respect to $(t,f)$.
\end{proposition}

\begin{proof}
Let $t\in\cI$ and $f$, $\delta_f \in Y$ such that $\delta_f \ne 0$.  Then
\begin{align*}
  F(t,f+\delta_f) &= \widetilde{F}[t,f+\delta_f,f+\delta_f]
  \\[1ex]
  &= \widetilde{F}[t,f,f] + \widetilde{F}[t,f,\delta_f] + \widetilde{F}[t,\delta_f,f]
  + \widetilde{F}[t,\delta_f,\delta_f]
  \\[1ex]
  &= F(t,f) + \widetilde{F}[t,f,\delta_f] + \widetilde{F}[t,\delta_f,f] + F(t,\delta_f).
\end{align*}
So equation \eqref{C42} is satisfied with
\begin{equation}\label{C47}
  \bigl[D_YF(t,f)\bigr]\delta_f = \widetilde{F}[t,f,\delta_f] + \widetilde{F}[t,\delta_f,f]
\end{equation}
and $R_Y\bigl((t,f),\delta_f\bigr)=F(t,\delta_f)$.

Relation \eqref{C30} implies that
\[
  \big\|\widetilde{F}[t,f,\delta_f]+\widetilde{F}[t,\delta_f,f]\big\|_X
  \le 2c(t)\|f\|_Y\|\delta_f\|_Y.
\]
Thus, for each fixed $(t,f)\in\cI\times Y$, we have that $D_YF(t,f)$ is a
bounded linear operator in $Y$.

Take $t'\in\cI$ and fix $f\in Y$.  Then, for $\delta_f\in Y\setminus\{0\}$
and $t\in[t_0,t']$, we have
\[
  \frac{\big\|R_Y\bigl((t,f),\delta_f\bigr)\big\|_X}{\|\delta_f\|_Y}
  = \frac{\|F(t,\delta_f)\|_X}{\|\delta_f\|_Y}
  = \frac{\|\widetilde{F}[t,\delta_f,\delta_f]\|_X}{\|\delta_f\|_Y}
  \le c(t')\|\delta_f\|_Y
  \to 0
\]
as $\|\delta_f\|_Y \to 0$, and the convergence is uniform in $t \in [t_0,t']$.
Thus $F$ is Fr\'echet differentiable with respect to the second argument,
and the Fr\'echet derivative at $f \in Y$ is given by \eqref{C47}
for all $\delta_f \in Y$.
Hence assumption (a) of Lemma~\ref{C41} holds.

Fix $(t,f)\in\mr{\cI}\times Y$.
From assumption~(a) of the current proposition we have, for $\delta_t\in\RR\setminus\{0\}$
such that $t+\delta_t\in\mr{\cI}$,
\begin{align*}
  F(t+\delta_t,f) &= \widetilde{F}[t+\delta_t,f,f]
  = \widetilde{F}[t,f,f] + \Bigl(\frac{\partial}{\partial t}\widetilde{F}[t,f,f]\Bigr)\delta_t 
  + R\bigl((t,f,f),\delta_t\bigr)
  \\[1ex]
  &= F(t,f) + \Bigl(\frac{\partial}{\partial t}\widetilde{F}[t,f,f]\Bigr)\delta_t 
  + R\bigl((t,f,f),\delta_t\bigr),
\end{align*}
where
\[
  \frac{\big\|R\bigl((t,f,f),\delta_t\bigr)\|_X}{\|\delta_t\|} \to 0
  \qquad \text{as} \;\; |\delta_t| \to 0.
\]
Thus $F$ is strongly differentiable with respect to the first argument,
and so assumption (b) in Lemma~\ref{C41} holds.

From Lemma~\ref{C29}\,(ii) we have
that $\widetilde{F}$ is continuous.  It follows from \eqref{C47}
that $D_YF(t,f)$ is also continuous with respect to $(t,f)$.
Using assumption (b) and the continuity of the derivative from (a)
we can prove in the same way as in the second half of the proof of
Lemma~\ref{C29}
that $(t,f,g) \mapsto \frac{\partial}{\partial t} \widetilde{F}[t,f,g]$ is continuous,
and hence $(t,f) \mapsto \frac{\partial}{\partial t} F(t,f)$ is continuous,
i.e.\ (c) in Lemma~\ref{C41} is satisfied.

Finally, assumption (d) in Lemma~\ref{C41}
follows from \eqref{C47} and the current assumption (b).
Now the assertions follow from Lemma~\ref{C41}.
\end{proof}

\subsection{Positive semigroups}
\label{C12}

Let us recall some definitions from the theory of positive $C_0$-semigroups,
which will be crucial for our investigations.
For a detailed theory of general $C_0$-semigroups and positive semigroups, in particular,
we refer the reader to \cite{engel_nagel2000,banasiak_arlotti2006perturbations,aliprantis_burkinshaw2006,batkai_kramer_rhandi2017positive}.

Let $X$ be a real \emph{ordered vector space}. 
The set of non-negative elements in $X$ is known as the \emph{positive cone} and
is denoted by $X_+$.  Likewise, for a subspace $D$ of $X$, the set of
non-negative elements in $D$ is denoted by $D_+$.
Further, a linear mapping $A:Y\to X$, where $X,Y$ are ordered vector spaces, 
is said to be \emph{positive} if $Ay\ge0$ for all $y\in Y_+$.

Now suppose that $X$ is a \emph{vector lattice} (or \emph{Riesz space}), i.e.\ it is an ordered vector space 
and $\sup\{f,g\}$ exists in $X$ for all $f,g\in X$.  
Let $f \in X$.
Then $f_{\pm} \coloneqq \sup\{\pm f,0\}$ and $|f| \coloneqq \sup\{f,-f\}$ are well defined 
and satisfy $f=f_+-f_-$ and $|f|=f_++f_-$.
A vector lattice is a \emph{Banach lattice} if if is a Banach space with norm $\normcdot$
and $|x|\le|y|$ implies $\|x\|\le\|y\|$ for all $x,y\in X$.
Note that in a Banach lattice the relation $\||x|\|=\|x\|$ is valid for all $x\in X$;
see \cite[(2.74)]{banasiak_arlotti2006perturbations}.

If $X$ is a Banach lattice and satisfies $\|f+g\|=\|f\|+\|g\|$
for all $f,g\in X_+$, then $X$ is called an \emph{AL-space}.
It follows from \cite[Theorems~2.64 and 2.65]{banasiak_arlotti2006perturbations} that, 
when $X$ is an AL-space, there is a unique bounded linear functional, $\phi$, 
that extends $\normcdot$ from $X_+$ to~$X$.

In the following lemma we consider two positive operators where one dominates
the other with respect to the order structure;
this is a slight generalisation of \cite[Theorem~2.68]{banasiak_arlotti2006perturbations}.

\begin{lemma}\label{C48}
Let $X$, $Y$ and $Z$ be Banach lattices such that $Y$ is a sublattice of $X$.
Moreover, let $A:Y\to Z$ and $B:X\to Z$ be positive operators such that
\[
  Ay \le By \qquad \text{for all} \;\; y\in Y_+.
\]
Then $A$ is bounded as an operator from $(Y,\normcdotsub{X})$ to $(Z,\normcdotsub{Z})$,
and $B\in\cB(X,Z)$.
If $Y$ is dense in $X$ and $\overline{A}$ denotes the closure of $A$ as an operator
from $X$ to $Z$, then
\[
  \big\|\overline{A}\big\|_{\cB(X,Z)} \le \big\|B\big\|_{\cB(X,Z)}.
\]
\end{lemma}

\begin{proof}
It follows from \cite[Theorem~2.65]{banasiak_arlotti2006perturbations} that $B$ is bounded.
Let $y\in Y$.  Since $Y$ is a sublattice of $X$, the moduli in $X$ and $Y$ coincide
on $Y$, and therefore $|y|\in Y_+$.
The positivity of $A$ and the relation $-|y|\le y\le|y|$ imply that
$-A|y|\le Ay\le A|y|$ and hence
\[
  |Ay| = \sup\{Ay,-Ay\} \le A|y| \le B|y| = \big|B|y|\big|.
\]
From this we can deduce that
\[
  \|Ay\|_Z \le \big\|B|y|\big\|_Z \le \|B\|_{\cB(X,Z)}\big\||y|\big\|_X
  = \|B\|_{\cB(X,Z)}\|y\|_X,
\]
and all the assertions follow.
\end{proof}

Let $(S(t))_{t \ge 0}$ be a $C_0$-semigroup on an ordered Banach space $X$.
Then $(S(t))_{t \ge 0}$ is a \emph{positive semigroup} if $S(t)f \ge 0$
for all $f \ge 0$, $t \ge 0$.
A positive semigroup $(S(t))_{t\ge0}$ is called \emph{substochastic} if $\|S(t)\|\le 1$ 
for all $t\ge0$, and \emph{stochastic} if $\|S(t)f\|=\|f\|$ for all $t\ge0$ and $f\in X_+$.

\subsection{Analytic semigroups and interpolation spaces}
\label{C13}

Analytic semigroups, see \cite[Definition~II.4.5]{engel_nagel2000},
are of particular importance in Section~\ref{C04} to weaken assumptions on the coagulation kernel.
These semigroups have very useful properties that make them desirable to work with.
In particular, they allow us to work in interpolation spaces.
Let us recall their definition.

Let $X$ be a Banach space and let $G$ be a densely defined operator
such that $G$ is the generator of an analytic $C_0$-semigroup, $(S(t))_{t\ge0}$, on $X$
(for this one has to use a complexification of the space; 
for details see, e.g.\ \cite[Section~2.2.5]{banasiak_arlotti2006perturbations}).
Consider $\cD(G)$ equipped with the graph norm.

A Banach space $(Y,\normcdotsub{Y})$ with the property
$\cD(G) \hookrightarrow Y \hookrightarrow X$,
where $\hookrightarrow$ denotes continuous embedding,
is called an \emph{intermediate space} between $\cD(G)$ and $X$.
If, in addition, for all $T\in\cB(X)$ such that $T|_{\cD(G)}\in\cB(\cD(G))$,
we have $T|_{Y}\in\cB(Y)$, then $Y$ is called an \emph{interpolation space}
between $\cD(G)$ and $X$.

As in \cite[\S 2.2.1]{lunardi1995analytic}, for $\alpha\in(0,1)$ and $p\in[1,\infty]$
we define a class of intermediate spaces between $X$ and $\cD(G)$ by
\begin{align}
  D_G(\alpha,p) &\coloneq \Bigl\{x \in X:
  t \mapsto v(t)\coloneq\big\|t^{1-\alpha-\frac{1}{p}}GS(t)x\big\|_X \in L^p(0,1)\Bigr\},
  \label{C49}
  \\[1ex]
  \|x\|_{D_G(\alpha,p)} &\coloneq \|x\|_X + \|v\|_{L^p(0,1)},
  \label{C50}
\end{align}
where $v$ in \eqref{C50} is as
in \eqref{C49}.  Further, we define
\[
  D_G(\alpha) \coloneqq \bigl\{x\in D_G(\alpha,\infty):
  \lim_{t\to0}t^{1-\alpha}GS(t)x=0\bigr\},
\]
and, for all $p\in[1,\infty]$, we set
\begin{alignat*}{2}
  D_G(0,p) &\coloneqq X, \qquad & \normcdotsub{D_G(0,p)} &\coloneqq \normcdotsub{X},
  \\[1ex]
  D_G(0) &\coloneqq X, \qquad & \normcdotsub{D_G(0)} &\coloneqq \normcdotsub{X}.
\end{alignat*}
It follows from \cite[Propositions~2.2.2 and 1.2.6]{lunardi1995analytic}
that $D_G(\alpha,p)$ are (real) interpolation spaces.
We note that, by \cite[Proposition~2.2.7]{lunardi1995analytic}, the part of $G$ in $D_G(\alpha,p)$
is sectorial in $D_G(\alpha,p)$ and so is also the generator of
an analytic $C_0$-semigroup on $D_G(\alpha,p)$.

We also need the notion of fractional powers of $-G$ under the extra assumption that
the semigroup generated by $G$ has a strictly negative growth bound;
see, e.g.\ \cite[\S 2.2.2]{lunardi1995analytic}.
For $\alpha>0$ we define negative fractional powers by
\[
  (-G)^{-\alpha} \coloneq \frac{1}{\Gamma(\alpha)}
  \int\limits_0^{\infty} t^{\alpha-1}S(t)\,\rd t,
\]
and positive powers by
\begin{equation}\label{C51}
  (-G)^{\alpha} \coloneq \bigl((-G)^{-\alpha}\bigr)^{-1}.
\end{equation}
Moreover, we set $G^0\coloneq I$.

In the following proposition some relations among these interpolations spaces are summarised.

\begin{proposition}\label{C52}
Let $G$ be the generator of an analytic semigroup on the space $X$.
\begin{myenum}
\item
For $\alpha\in(0,1)$ and $p_1,p_2\in[1,\infty)$ with $p_1<p_2$ we have
\[
	 D_G(\alpha,p_1) \hookrightarrow D_G(\alpha,p_2) \hookrightarrow D_G(\alpha) \hookrightarrow D_G(\alpha,\infty).
\]
\item
For $\alpha_1,\alpha_2\in(0,1)$ with $\alpha_1<\alpha_2$ we have
\[
	D_G(\alpha_2,\infty) \hookrightarrow D_G(\alpha_1,1).
\]
\item
Assume, in addition, that the semigroup generated by $G$ has a strictly negative growth bound.
Then, for $\alpha\in(0,1)$,
\[
	D_G(\alpha,1) \hookrightarrow \cD\bigl((-G)^{\alpha}\bigr) \hookrightarrow D_G(\alpha,\infty).
\]
\end{myenum}
\end{proposition}

\begin{proof}
Items (i) and (iii) follow from \cite[Corollary~2.2.3]{lunardi1995analytic}
and \cite[Proposition~2.2.15]{lunardi1995analytic} respectively.
The statement in (ii) follows from \cite[Proposition~1.2.3]{lunardi1995analytic} combined with
\cite[Proposition~2.2.2]{lunardi1995analytic}.
\end{proof}

\medskip

\noindent
Finally, let us recall the definition of spaces of H\"older continuous functions.
Let $X$ be a Banach space.  For each $\gamma\in(0,1)$ and an interval $\cI \subseteq \RR$,
let $C_b(\cI,X)$ be the set of bounded, continuous functions from $\cI$ into $X$,
and define the spaces of H\"older continuous functions by
\begin{equation}\label{C95}
  C^{\gamma}(\cI,X) \coloneqq \biggl\{f \in C_b(\cI,X):
  [f]_{C^{\gamma}(\cI,X)}\coloneqq\sup\limits_{t,s\in \cI: \, s<t}
  \frac{\|f(t)-f(s)\|}{(t-s)^{\gamma}}<\infty\biggr\}
\end{equation}
with norm
\begin{equation*}
  \|f\|_{C^{\gamma}(\cI,X)} \coloneqq \|f\|_{C_b(\cI,X)} + [f]_{C^{\gamma}(\cI,X)}.
\end{equation*}

\section{Abstract Existence, Uniqueness and Positivity \\ Results}\label{C03}

In \cref{C31,C32}, the main results of this section, we prove existence, uniqueness and 
positivity of mild solutions of semi-linear, abstract Cauchy problems in a Banach space or Banach lattice $X$.
One main novel feature is the use of a subspace $Y$ of $X$, which can be
a subspace of elements of some sort of higher regularity.
We allow the non-linearity in the equation to be more singular in the sense
that it is continuous only from $Y$ to $X$.
This is compensated by an assumption that the semigroup connected with the 
linear part of the equation is regularity enhancing.
Another important result, presented in \cref{C32}, is the positivity of solutions 
in the setting of Banach lattices under relatively weak assumptions.

Let $X$, $Y$ be Banach spaces such that $Y$ is continuously embedded in $X$.
We consider the following semi-linear, abstract Cauchy problem~(ACP)
\begin{align}
  u'(t) &= Gu(t)+F\bigl(t,u(t)\bigr), \quad t \in (t_0,T), 
  \label{C54}
  \\[1ex]
  u(t_0) &= \mathring{u},
  \label{C55}
\end{align}
where
\begin{mybull}
\item
	$t_0<T \le \infty$;
\item
	$G$ is the generator of a $C_0$-semigroup, $(S(t))_{t\ge0}$, on $X$;
\item
	$F: [t_0, T) \times Y \to X$ is continuous;
\item
	$\mr{u}\in Y$.
\end{mybull} 

\noindent
We recall the following standard notions of solutions of such abstract Cauchy problems;
see, e.g.\ \cite[Definitions~4.2.1 and 6.1.1]{pazy1983semigroups}.

\begin{definition}\label{C56}
Let $X$, $Y$, $G$, $(S(t))_{t\ge0}$, $F$ and $\mr{u}$ be as above.
\begin{myenum}
\item
	A function $u\in C([t_0,t_1),Y)$ with $t_1\in(t_0,T]$ is said to be a \emph{mild solution}
	of \eqref{C54}, \eqref{C55} if $u$ satisfies the integral equation
	\begin{equation}\label{C57}
		u(t) = S(t-t_0)\mathring{u}+\int\limits_{t_0}^t S(t-s)F(s,u(s))\,\rd s,
		\qquad  t \in [t_0,t_1).
	\end{equation}
\item
	We say that $u$ is a \emph{classical solution} of \eqref{C54}, \eqref{C55} 
	on $[t_0,t_1)$ if $u:[t_0,t_1) \to Y$ is continuous,
	$u:(t_0,t_1) \to X$ is continuously differentiable, 
	$u(t) \in \cD(G)$ for $t \in (t_0,t_1)$ and \eqref{C54} and \eqref{C55} are satisfied.
\item
	A mild (respectively classical) solution, $u$, of \eqref{C54}, \eqref{C55} 
	on $[t_0,t_1)$ is the \emph{maximal mild solution} (respectively \emph{maximal classical solution})
	if there does not exist a $\tilde{t}_1\in(t_1,T]$
	and an extension $\tilde{u}$ of $u$ such that $\tilde{u}$ is a mild
	(respectively classical) solution of \eqref{C54}, \eqref{C55}
	on $[t_0,\tilde{t}_1)$.
\end{myenum}
\end{definition}

\subsection{General existence, uniqueness and positivity results}
\label{C14}

We formulate the following main theorems of this section, \cref{C31,C32}, which are applied to the
coagulation--fragmentation system in Section~\ref{C04}.
These theorems are based on ideas that go back to \cite{fujita_kato1964} and \cite[Section~3.3]{henry1981} and are also covered in
\cite[Theorem~6.1.4]{pazy1983semigroups} and \cite[Theorem~7.1.2]{lunardi1995analytic}, and we use a similar method to prove
these results.
However, unlike previous theorems, we allow the non-linear operator $F$
in \eqref{C54} to map from an arbitrary Banach space $Y$
into another Banach space, $X$, with $Y$ continuously embedded in $X$.
The case when $Y$ is an interpolation space is discussed in more detail in \cref{C86}.
Moreover, we obtain a positivity result in Theorem~\ref{C32}
that is absent in \cite[Theorem~6.1.4]{pazy1983semigroups} 
and \cite[Theorem~7.1.2]{lunardi1995analytic}.
This positivity result is also one reason for presenting the proof of \cref{C31}
because we have to track the positivity of various mappings that arise 
in the proof of \cref{C31}.  We begin in this subsection, by presenting the statements of
Theorem~\ref{C31} and \ref{C32}.  We then prove the former in Section~\ref{C16}, following the presentation
of two necessary lemmas.  In Section~\ref{C17}, we present lemmas required in the proof of 
Theorem~\ref{C32}, before providing the proof of this theorem itself.  Finally, in Section~\ref{C15} we consider
some special situations and additions.

In the following theorem, the first main result of this section, existence and uniqueness of 
solutions of abstract semi-linear ACPs are established.

\begin{theorem}\label{C31}
Let $(X,\normcdotsub{X})$, $(Y,\normcdotsub{Y})$ be Banach spaces such that $Y$
is continuously embedded in $X$ and let $T\in(0,\infty]$.  Assume that
\begin{myenuma}
\item
	$G$ is the generator of a $C_0$-semigroup, $(S(t))_{t\ge0}$, on $X$ such that $S(t)$
	leaves $Y$ invariant for every $t \in [0,\infty)$,
	and that $t \mapsto S(t)y$ is continuous from $[0,\infty)$ to $Y$ for every $y \in Y$;
\item
	there exists a function $\eta: (0,\infty) \to (0,\infty)$
	such that $\eta(\delta) \to 0$ as $\delta \to 0^+$ and, for all $s_0,s_1\in[0,\infty)$
	with $s_0<s_1$ and $\varphi \in C([s_0,s_1],X)$, we have
	\begin{align}
		& \int\limits_{s_0}^{s_1} S(s_1-s)\varphi(s)\,\rd s \in Y,
		\label{C58}
		\\[1ex]
		& \Bigg\|\int\limits_{s_0}^{s_1} S(s_1-s)\varphi(s)\,\rd s \Bigg\|_Y
		\le \eta(s_1-s_0)\big\|\varphi\big\|_{C([s_0,s_1],X)}\,;
		\label{C59}
	\end{align}
\item
	$F:[0,T) \times Y \to X$ is continuous 
	and satisfies a Lipschitz condition in the second argument on bounded sets,
	uniformly in the first argument on compact intervals.
\end{myenuma}
Let $\mr u\in Y$.  Then the following statements hold.
\begin{myenum}
\item
	There exists a unique maximal mild solution $u\in C([0,\tmax),Y)$
	of \eqref{C54}, \eqref{C55} \textup{(}with $t_0=0$\textup{)}
	where $\tmax\in(0,T]$.
\item
	If $\tmax<T$, then $\|u(t)\|_Y \to \infty$ as $t\to\tmax^-$.
\item
	For every $\tau\in(0,\tmax)$ there exist $r_\tau>0$ and $M_\tau>0$ such that,
	whenever $\mr v\in\overline{B}_Y(\mr u,r_\tau)$ and $v$ is the corresponding
	maximal mild solution, then $v$ is defined on $[0,\tau]$ and
	\begin{equation}\label{C60}
		\|v-u\|_{C([0,\tau],Y)} \le M_\tau\|\mr v-\mr u\|_Y\,.
	\end{equation}
\end{myenum}
\end{theorem}

\begin{remark}\label{C61}
\rule{0ex}{1ex}
\begin{myenum}
\item
	Note that assumptions (a) and (b) in \cref{C31} concern the operator $G$ 
	in \eqref{C54} whereas assumption (c) is about the non-linear term $F$.
\item
	The expression in \eqref{C58}, as a function of $s_1$,
	is a mild solution of the linear inhomogeneous ACP $u'(t)=Gu(t)+\varphi(t)$, $u(s_0)=0$.
\item
	That an inequality of the form \eqref{C59} holds, follows from \eqref{C58}
	and the Closed Graph Theorem.  What is assumed additionally in (b) is
	that $\eta(\delta)\to0$ as $\delta\to0$.
\item
	The local Lipschitz property in item (iii) seems to be new in this form.
	From this one can easily derive a statement about convergent sequences of initial conditions.
\item\label{C130}
	Assume that $Y=X$ and that $G$ generates a $C_0$-semigroup on $X$.
	Then (a) and (b) are satisfied as can be easily seen.
\item
	In \cref{C86} below it is shown that (a) and (b) are satisfied in the case when $Y$
	is a certain interpolation space.
\item
	One can use \cref{C29} to check assumption (c).
\end{myenum}
\end{remark}

\medskip

\noindent
The following theorem establishes the above mentioned positivity of the mild solution in Theorem~\ref{C31}\,(i)
under relatively weak assumptions.  The proof of this Theorem, which we provide in \S \ref{C17}, 
relies on the addition and subtraction of a positive operator, $H$, to the C--F system.  We believe that this idea originates 
from work in a different context by Belleni-Morante \cite[Section~8.6]{belleni1994concise}, 
and it has been used in prior investigations of the C--F system, e.g.\ for uniformly bounded coagulation 
in \cite{mcbride_smith_lamb2010strongly} and for unbounded coagulation in \cite{banasiak_lamb2012analytic}.

\begin{theorem}\label{C32}
Let assumptions \textup{(a)--(c)} of Theorem~\ref{C31}
hold.  In addition, assume that the following conditions are satisfied:
\begin{myenuma}
\setcounter{counter_b}{3}
\item
  $X$ and $Y$ are Banach lattices such that $Y$ is a sublattice of $X$ and
  that $\cD(G)$ is a dense subspace of $Y$;
\item
  there exists a positive linear operator $H:Y\to X$ such that,
  for all $\gamma\ge0$, $G-\gamma H$ generates a positive semigroup
  $(S_{\gamma}(t))_{t\ge0}$ on $X$ and $S_\gamma(t)$ leaves $Y$ invariant
  for every $t\in[0,\infty)$;
\item
  for every $\tau\in[0,T)$ and every $r>0$ there exists $\gamma\ge0$ such that
  \begin{equation}\label{C62}
    F(t,v)+\gamma Hv \ge 0 \qquad \text{for all} \;\; v \in \overline{B}_Y(0,r)_+, \
    t \in [0,\tau].
  \end{equation}
\end{myenuma}
Let $\mr{u}\in Y_+$ and let $u$ be the maximal mild solution
on $[0,\tmax)$ from \cref{C31} with $\tmax\in(0,T]$.  Then $u(t)\ge0$ for all $t\in[0,\tmax)$.
\end{theorem}

\begin{remark}\label{C97}
\rule{0ex}{1ex}
\begin{myenum}
\item
	In Proposition~\ref{C86} and Lemma~\ref{C87} below we provide some sufficient conditions for 
	assumptions \textup{(d)} and \textup{(e)} to hold.
\item
	If $F(t,\cdot)$ maps positive elements onto positive elements, then one can choose $H=0$,
	in which case assumption (e) simplifies to the assumption that $(S(t))_{t\ge0}$ is a positive semigroup.
\item
	Assume that $X$ is a Banach lattice, $Y=X$ and $G$ generates a positive semigroup.
	Then (d) and (e) are satisfied with $H=I$ since then $S_\gamma(t)=e^{-\gamma t}S(t)$.
\end{myenum}
\end{remark}

\subsection{Proof of \cref{C31}}\label{C16}

Before we prove \cref{C31}, we need some lemmas.

\begin{lemma}\label{C63}
Let $X$ and $Y$ be Banach spaces such that $Y$ is continuously embedded in $X$.
Further, let $G$ be the generator of a $C_0$-semigroup, $(S(t))_{t\ge0}$, on $X$ 
such that $S(t)$ leaves $Y$ invariant for every $t\in[0,\infty)$.
\begin{myenum}
\item
  For every $t\in[0,\infty)$ we have $S(t)\restrict_Y\in\cB(Y)$.
\item
  If $t\mapsto S(t)y$ is continuous from $[0,\infty)$ to $Y$ for every $y\in Y$, then
  \begin{equation}\label{C64}
    \bigl\{\big\|S(t)\restrict_Y\big\|_{\cB(Y)}: t \in [0,t_1]\bigr\}
  \end{equation}
  is bounded for every $t_1\in(0,\infty)$.
\item
  Assume that the set in \eqref{C64} is bounded for every $t_1\in(0,\infty)$ and 
  that assumption \textup{(b)} in Theorem~\ref{C31} holds.
  Then, for $s_0,s_1\in[0,\infty)$ with $s_0<s_1$ and every $\varphi \in C([s_0,s_1],X)$,
  the function
  \[
    t \mapsto \int\limits_{s_0}^t S(t-s)\varphi(s)\,\rd s
  \]
  is continuous from $[s_0,s_1]$ to $Y$.
\end{myenum}
\end{lemma}

 \begin{proof}
(i)
Let $t\in[0,\infty)$.  Since $S(t)$ is closed in $X$, the space $Y$ is continuously 
embedded in $X$ and $Y$ is invariant under $S(t)$, it follows that the 
operator $S(t)\restrict_Y$ is closed in $Y$.
Hence the Closed Graph Theorem implies that $S(t)\restrict_Y$ is bounded.

(ii)
Let $t_1\in(0,\infty)$.  For every $y\in Y$, the mapping $t \mapsto S(t)y$ is continuous 
from $[0,t_1]$ into $Y$  and hence $\{\|S(t)y\|: t\in[0,t_1]\}$ is bounded.
The Uniform Boundedness Principle implies that the set in \eqref{C64} is bounded.

(iii)
Without loss of generality we can assume that the function $\eta$ is increasing
on $(0,\infty)$.
Let $s_0,s_1\in[0,\infty)$ with $s_0<s_1$ and let $\varphi\in C([s_0,s_1],X)$.
Further, let $t_1,t_2\in[s_0,s_1]$ such that $t_1<t_2$,
and set $h\coloneqq t_2-t_1$.  Then
\begin{align*}
  & \left\|\,\int\limits_{s_0}^{t_2} S(t_2-s)\varphi(s)\,\rd s
  - \int\limits_{s_0}^{t_1} S(t_1-s)\varphi(s)\,\rd s\,\right\|_Y
  \\[1ex]
  &= \left\|\,\int\limits_{s_0}^{t_2} S(t_2-s)\varphi(s)\,\rd s
  - \int\limits_{s_0+h}^{t_2} S(t_2-s)\varphi(s-h)\,\rd s\,\right\|_Y
  \displaybreak[0]\\[1ex]
  &\le \left\|\,\int\limits_{s_0}^{s_0+h} S(t_2-s)\varphi(s)\,\rd s\,\right\|_Y
  + \left\|\;\int\limits_{s_0+h}^{t_2} S(t_2-s)
  \bigl[\varphi(s)-\varphi(s-h)\bigr]\,\rd s\,\right\|_Y
  \displaybreak[0]\\[1ex]
  &\le \big\|S(t_2-s_0-h)\big\|_{\cB(Y)}\,
  \left\|\,\int\limits_{s_0}^{s_0+h} S(s_0+h-s)\varphi(s)\,\rd s\,\right\|_Y
  \\
  &\quad + \eta(t_2-s_0-h)\big\|\varphi(\cdot)-\varphi(\cdot-h)\big\|_{C([s_0+h,t_2],X)}
  \\[1ex]
  &\le \sup_{t \in [0,s_1]}\big\|S(t)\restrict_Y\big\|_{\cB(Y)}\,
  \eta(h)\|\varphi\|_{C([s_0,s_1],X)}
  + \eta(s_1-s_0)\big\|\varphi(\cdot)-\varphi(\cdot-h)\big\|_{C([s_0+h,s_1],X)}.
\end{align*}
The first term clearly converges to $0$ as $h\to0$.
For the second term note that $\varphi$ is uniformly
continuous on $[s_0,s_1]$ as it is continuous on a compact interval.
Hence $\|\varphi(\cdot)-\varphi(\cdot-h)\big\|_{C([s_0+h,s_1],X)}\to0$
as $h\to0$.
\end{proof}

We now introduce some notation that is used in the following.
Let assumptions (a)--(c) of Theorem~\ref{C31} hold.
We assume again, without loss of generality, that the function $\eta$ is increasing
on $(0,\infty)$.
Further, we define
\begin{align}
  \hatM(t) &\coloneqq \sup\Bigl\{\big\|S(s)\restrict_Y\big\|_{\cB(Y)}:
  s \in [0,t]\Bigr\}, \qquad t\in(0,T),
  \label{C65}
  \\[1ex]
  \hatM_0 &\coloneqq \limsup\limits_{t \to 0^+} \hatM(t).
  \label{C66}
\end{align}
By Lemma~\ref{C63}\,(ii), $\hatM(t)$ is well defined
for each $t\in(0,T)$.  Since $S(0)\restrict_Y=I\restrict_Y$,
we have $\hatM(t)\ge1$ for all $t\in(0,T)$ and hence $\hatM_0\ge1$.
Moreover, for $t \in (0,T)$, we set
\begin{equation}\label{C67}
  N(t) \coloneqq \max\bigl\{\|F(s,0)\|_X:s\in[0,t]\bigr\}.
\end{equation}
Then $N$ and $\hatM$ are both monotone increasing functions.
Further, let $L(\tau,r)$ be the Lipschitz constant associated with $F$
as in Definition~\ref{C27}.
Without loss of generality, we assume that $L$ is monotone increasing 
in the first argument.

For $C>0$, $\tau\in(0,T)$ and $r>2C\hatM_0$, we define
\begin{equation}\label{C68}
  \delta(C,r,\tau) \coloneqq \frac{1}{2}\sup\biggl\{d\in(0,\tau):
  L(\tau,r)\eta(d) \le \frac{1}{2} \;\;\text{and}\;\;
  C\hatM(d)+N(\tau)\eta(d) \le \frac{r}{2}\biggr\}.
\end{equation}
The set on the right-hand side of \eqref{C68} is non-empty because
\[
  \lim_{d\to0^+}L(\tau,r)\eta(d) = 0
  \qquad\text{and}\qquad
  \limsup_{d\to0^+}\Bigl(C\hatM(d)+N(\tau)\eta(d)\Bigr)
  = C\hatM_0 < \frac{r}{2}\,,
\]
and hence $\delta=\delta(C,r,\tau)>0$.
Moreover, the factor of $\frac{1}{2}$ in front of $\sup$ in \eqref{C68} and the
monotonicity of $\eta$ and $\hatM$ imply that $\delta$ is contained in the set on the
right-hand side of \eqref{C68}.  It follows that
\begin{equation}\label{C69}
  L(\tau,r)\eta(\delta) \le \frac{1}{2} \qquad \text{and} \qquad
  C\hatM(\delta)+N(\tau)\eta(\delta) \le \frac{r}{2},
\end{equation}
where again $\delta=\delta(C,r,\tau)$.
Let $t_0,t_1\in(0,T)$ be such that $t_0<t_1$.
Consider the Banach space 
\[
  Z \coloneqq C([t_0,t_1],Y) \qquad\text{with norm} \ 
  \|v\|_Z=\max\limits_{t\in[t_0,t_1]}\|v(t)\|_Y, \;\; v \in Z,
\]
and the set
\begin{equation}\label{C70}
  \Sigma_r \coloneqq \overline{B}_Z(0,r)
  = \bigl\{v \in Z: \|v\|_Z \le r\bigr\}.
\end{equation}
For $\mr u\in Y$, define the mapping $Q_{\mr u}$ by
\begin{equation}\label{C71}
  (Q_{\mr u}v)(t) \coloneqq S(t-t_0)\mathring{u}+\int\limits_{t_0}^t
  S(t-s)F\bigl(s,v(s)\bigr)\,\rd s,
  \qquad v \in \Sigma_r, \ t_0 \le t \le t_1.
\end{equation}
We use this notation in the following key lemma.
The proof of part (i) of this lemma is based on ideas used in the proof
of \cite[Theorem~6.1.4]{pazy1983semigroups}.  However, the latter deals only with
the particular case where $Y=X$.
Moreover, we provide a positivity result that is absent in \cite[Theorem~6.1.4]{pazy1983semigroups}.

\begin{lemma}\label{C72}
Suppose that the assumptions \textup{(a)--(c)}
of Theorem~\ref{C31} are satisfied.
Let $C > 0$ and $r > 2C\hatM_0$, let $t_0,\tau\in[0,T)$ with $t_0<\tau$ 
and set
\begin{equation}\label{C73}
  t_1 = \min\bigl\{\tau,t_0+\delta(C,r,\tau)\bigr\},
\end{equation}
where $\delta(C,r,\tau)$ is as in \eqref{C68}.
Further, let $\Sigma_r$ and $Q_{\mr u}$ be defined as in \eqref{C70}
and \eqref{C71} respectively.

\begin{myenum}
\item
  For every $\mr u\in Y$ with $\|\mr u\|_Y\le C$ the mapping $Q_{\mr u}$
  is a contraction on $\Sigma_r$ with contraction constant at most $\frac{1}{2}$;
  hence it has a unique fixed point, $u$, in $\Sigma_r$.
\item
  Let $\mr u_1,\mr u_2\in Y$ with $\|\mr u_1\|_Y,\|\mr u_2\|_Y\le C$ and let
  $u_i$ be the corresponding fixed point of $Q_{\mr u_i}$ for $i\in\{1,2\}$.
  Then
  \begin{equation}\label{C74}
    \|u_1-u_2\|_Z \le 2\hatM(t_1-t_0)\|\mr u_1-\mr u_2\|_Y.
  \end{equation}
\item
  Suppose that, in addition, $X$ is an ordered Banach space, $(S(t))_{t\ge0}$
  is a positive semigroup and $\mr{u}\in Y_+$.  Assume that
  \[
    F(t,v) \ge 0 \qquad \text{for all} \;\; v \in \overline{B}_Y(0,r)_+, \;\;
    t \in [t_0,t_1].
  \]
  Let $u$ be the unique fixed point of $Q_{\mr u}$ from \textup{(i)}.  Then $u\in(\Sigma_r)_+$.
\end{myenum}
\end{lemma}

\begin{proof}
(i)
Since $C<r$ (recall that $\hatM_0\ge1$), we have $\mr{u}\in\Sigma_r$,
where $\mr{u}$ is considered as a constant function.
It follows from assumption (c) of Theorem~\ref{C31}
that $F:[0,T)\times Y\to X$ is continuous.
Hence Lemma~\ref{C63}\,(iii) implies that, for $v\in\Sigma_r$,
the function
\[
  t \mapsto (Q_0v)(t) = \int\limits_{t_0}^t S(t-s)F\bigl(s,v(s)\bigr)\,\rd s
\]
is continuous from $[t_0,t_1]$ into $Y$.  
Since $t \mapsto S(t-t_0)\mathring{u}$ is continuous from $[t_0,t_1]$ into $Y$
by assumption (a) of Theorem~\ref{C31},
we obtain $Q_{\mr u}v\in Z$ for all $v\in\Sigma_r$.

Let $v,w\in\Sigma_r$.  Using \eqref{C59} and \eqref{C69} we obtain
\begin{align}
  \|Q_{\mr u}v-Q_{\mr u}w\|_Z &= \max_{t\in[t_0,t_1]}\left\|\,\int\limits_{t_0}^t
  S(t-s)\Bigl[F\bigl(s,v(s)\bigr)-F\bigl(s,w(s)\bigr)\Bigr]\,\rd s\,\right\|_Y
  \nonumber\\[1ex]
  &\le \max_{t\in[t_0,t_1]}\biggl(\eta(t-t_0)
  \max_{s\in[t_0,t]}\big\|F\bigl(s,v(s)\bigr)-F\bigl(s,w(s)\bigr)\big\|_X\biggr)
  \nonumber\displaybreak[0]\\[1ex]
  &\le \eta(t_1-t_0)L(\tau,r)\max_{s\in[t_0,t_1]}\big\|v(s)-w(s)\big\|_Y
  \nonumber\displaybreak[0]\\[1ex]
  &\le \eta(\delta)L(\tau,r)\|v-w\|_Z
  \nonumber\\[1ex]
  &\le \frac{1}{2}\|v-w\|_Z.
  \label{C75}
\end{align}

Now let $v\in\Sigma_r$.  Then \eqref{C75}
and \eqref{C69} imply that
\begin{align*}
  \|Q_{\mr u}v\|_Z &\le \|Q_{\mr u}v-Q_{\mr u}0\|_Z + \|Q_{\mr u}0\|_Z
  \\[1ex]
  &\le \frac{1}{2}\|v-0\|_Z + \max_{t\in[t_0,t_1]}\left[\|S(t-t_0)\mr{u}\|_Y
  + \left\|\,\int\limits_{t_0}^{t} S(t-s)F(s,0)\,\rd s \,\right\|_Y\right]
  \displaybreak[0]\\[1ex]
  &\le \frac{r}{2} + \hatM(t_1-t_0)\|\mr{u}\|_Y + \eta(t_1-t_0)N(t_1)
  \\[1ex]
  &\le \frac{r}{2} + C\hatM(\delta) + \eta(\delta)N(\tau)
  \le r,
\end{align*}
which shows that $Q_{\mr u}$ maps $\Sigma_r$ into itself.
It follows from \eqref{C75} that $Q_{\mr u}$ is a contraction
on $\Sigma_r$ with contraction constant at most $\frac{1}{2}$.
Hence $Q_{\mr u}$ has a unique fixed point, $u$, in~$\Sigma_r$.

(ii)
Since $u_i(t)=S(t-t_0)\mr u_i+(Q_0u_i)(t)$ for $i\in\{1,2\}$, we have
\begin{align*}
  \|u_1-u_2\|_Z &\le \max_{t\in[t_0,t_1]}\big\|S(t-t_0)\restrict_{Y}\big\|_{\cB(Y)}
  \|\mr u_1-\mr u_2\|_Y + \|Q_0u_1-Q_0u_2\|_Z
  \\[1ex]
  &\le \hatM(t_1-t_0)\|\mr u_1-\mr u_2\|_Y + \frac{1}{2}\|u_1-u_2\|_Z,
\end{align*}
which proves \eqref{C74}.

(iii)
Since, by assumption, $(S(t))_{t\ge0}$ is a positive semigroup, $\mr u\in Y_+$
and $F(t,v(t))\ge0$ for $v\in(\Sigma_r)_+$ and $t\in[t_0,t_1]$,
it follows that $Q_{\mr u}$ maps $(\Sigma_r)_+$ into itself.
The relation $u=\lim\limits_{n\to\infty}Q_{\mr u}^n \mr u$ implies
that $u\in(\Sigma_r)_+$.
\end{proof}

We can now use Lemma~\ref{C72}
to prove Theorem~\ref{C31}.
We adopt the notation introduced in the discussion preceding
Lemma~\ref{C72}.

\begin{proof}[Proof of Theorem~\ref{C31}]
(i)
We construct a maximal mild solution iteratively:
if a solution $u$ has been constructed on $[0,\hat t]$, then we apply
Lemma~\ref{C72}\,(i) with $C\ge\|u(\hat t)\|_Y$,
$r>2C\hatM_0$ and $\tau\in(\hat t,T)$ to obtain a solution on a larger interval;
for the initial step we use $\hat t=0$.
The uniqueness of a solution in $Y$ follows also from
Lemma~\ref{C72}\,(i):
if $u_1,u_2$ are two distinct solutions defined on $[0,\tmax^{(1)})$ and $[0,\tmax^{(2)})$
respectively, we can apply that lemma with
$t_0=\inf\{t\in[0,\min\{\tmax^{(1)},\tmax^{(2)}\}):u_1(t)\ne u_2(t)\}$
to obtain a contradiction.

(ii)
Assume that $\tmax<T$ and that $\|u(t)\|_Y \nrightarrow \infty$ as $t\to\tmax^-$.
Then there exist $s_n \in [0,\tmax)$ such that $(s_n)_{n=1}^{\infty}$
is increasing, $s_n\to\tmax$ as $n\to\infty$ and $\|u(s_n)\|_Y \le C$
for some $C>0$.  Choose $\tau\in(\tmax,T)$, $r>2C\hatM_0$ and
set $\delta \coloneqq \delta(C,r,\tau)$.  For some $n\in\NN$ we have $s_n>\tmax-\delta$.
Choose $t_1$ as in \eqref{C73} with $t_0=s_n$.  Then $t_1>\tmax$.
From Lemma~\ref{C72}\,(i)
it follows that the solution can be extended beyond $\tmax$,
which contradicts the definition of $\tmax$.

(iii)
Let $\tau\in(0,\tmax)$.  Choose $C\ge2\max_{t\in[0,\tau]}\|u(t)\|_Y$ and $r>2C\hatM_0$.
Further, set $\delta=\delta(C,r,\tau)$, choose $n\in\NN$ with $n\ge\frac{\tau}{\delta}$
and set
\[
  s_j \coloneqq \frac{j\tau}{n}\,, \qquad j=0,\ldots,n.
\]
Then $s_{j+1}-s_j=\frac{\tau}{n}\le\delta$.
Moreover, we set
\[
  K \coloneqq 2\hatM\Bigl(\frac{\tau}{n}\Bigr), \qquad 
  M_\tau \coloneqq K^n, \qquad
  r_\tau \coloneqq \frac{C}{2K^n}.
\]
Let $\mr v\in\overline{B}_Y(\mr u,r_\tau)$ and let $v$ be the corresponding maximal mild solution.
We prove by induction that, for $j\in\{0,\ldots,n\}$, $v$ is defined at $s_j$ and
\[
  \|v(s_j)\|_Y \le C, \qquad \|v(s_j)-u(s_j)\|_Y \le K^j\|\mr v-\mr u\|_Y.
\]
For $j=0$ this is clear since
\[
  \|v(s_0)\|_Y \le \|\mr v-\mr u\|_Y + \|\mr u\|_Y \le r_\tau + \frac{C}{2} \le C.
\]
Assume now that the statement is true for some $j\in\{0,\ldots,n-1\}$.
Lemma~\ref{C72} implies
that $v$ exists on $[s_j,s_{j+1}]$ and that, for $t\in[s_j,s_{j+1}]$,
\begin{equation}\label{C76}
  \|v(t)-u(t)\|_Y \le K\|v(s_j)-u(s_j)\|_Y \le K^{j+1}\|\mr v-\mr u\|_Y
  \le M_\tau\|\mr v-\mr u\|_Y.
\end{equation}
From this we obtain
\begin{align*}
  \|v(s_{j+1})\|_Y &\le \|v(s_{j+1})-u(s_{j+1})\|_Y + \|u(s_{j+1})\|_Y
  \\[1ex]
  &\le K^{j+1}\|\mr v-\mr u\|_Y + \frac{C}{2}
  \le K^n r_\tau + \frac{C}{2} = C.
\end{align*}
Hence the statement is true for all $j\in\{0,\ldots,n\}$ and therefore $v$ exists 
on $[0,s_n]=[0,\tau]$.
The inequality in \eqref{C60} follows from \eqref{C76}.
\end{proof}

\subsection{Proof of \cref{C32}}\label{C17}

To prove Theorem~\ref{C32}, we require the following three lemmas.
The first one is a generalisation of a standard perturbation result;
see, e.g.\ \cite[Corollary~III.1.7]{engel_nagel2000}.

\begin{lemma}\label{C77}
Assume that
\begin{myenuma}
\item
  $G$ and $(S(t))_{t\ge0}$ satisfy assumption \textup{(a)}
  of Theorem~\ref{C31};
\item
  $\cD(G)$ is a dense subspace of $Y$;
\item
  $\gamma\in\RR$ and $H\in\cB(Y,X)$ are such that $G-\gamma H$ generates
  a $C_0$-semigroup, $(S_\gamma(t))_{t\ge0}$, on $X$.
\end{myenuma}
Then
\begin{equation}\label{C78}
  S(t)y = S_{\gamma}(t)y + \int\limits_0^t S_{\gamma}(t-s)\gamma HS(s)y\,\rd s
\end{equation}
for every $t\ge0$ and $y\in Y$.
\end{lemma}

\begin{proof}
The relation in \eqref{C78} is trivially satisfied for $t=0$.
Let now $t>0$.
First take $y\in\cD(G)$ and consider the function
\[
  s \mapsto \xi_y(s) \coloneqq S_{\gamma}(t-s)S(s)y, \qquad s \in [0,t]
\]
as a function from $[0,t]$ to $X$.
Since $\cD(G)=\cD(G-\gamma H)$ is invariant under $(S(t))_{t\ge0}$,
we can use \cite[Lemma~B.16]{engel_nagel2000} to deduce that $\xi_y$ is differentiable
on $(0,t)$ and that
\begin{align*}
  \frac{\rd}{\rd s}\xi_y(s) &= -S_{\gamma}(t-s)(G-\gamma H)S(s)y + S_{\gamma}(t-s)GS(s)y
  \\[1ex]
  &= S_{\gamma}(t-s)\gamma HS(s)y.
\end{align*}
Integrating with respect to $s$ over the interval $[0,t]$ we can deduce that
\eqref{C78} is satisfied for $y\in\cD(G)$.

Now let $y\in Y$.  Since, by assumption, $\cD(G)$ is dense in $Y$,
there exist $y_n\in\cD(G)$, $n\in\NN$, such that $y_n\to y$ in $Y$.
Assumption (a) in Theorem~\ref{C31} implies
that $S(\cdot)y_n,S(\cdot)y\in C([0,t],Y)$,
and Lemma~\ref{C63}\,(ii) shows
that $S(\cdot)y_n\to S(\cdot)y$ in $C([0,t],Y)$.
Since $H\in\cB(Y,X)$, we therefore have $HS(\cdot)y_n\to HS(\cdot)y$
in $C([0,t],X)$.
From the first part of the proof we know that \eqref{C78} holds with $y$ replaced by $y_n$.
Taking the limit in $X$ on both sides we obtain \eqref{C78} also for $y$.
\end{proof}

In the next lemma the expression in \eqref{C57} is related to a similar
expression with $S$ replaced by $S_\gamma$.

\begin{lemma}\label{C79}
Assume that
\begin{myenuma}
\item
  $G$ and $(S(t))_{t\ge0}$ satisfy assumptions \textup{(a)} and \textup{(b)}
  of Theorem~\ref{C31};
\item
  $\cD(G)$ is a dense subspace of $Y$;
\item
  $\gamma\in\RR$ and $H\in\cB(Y,X)$ are such that $G-\gamma H$ generates
  a $C_0$-semigroup, $(S_\gamma(t))_{t\ge0}$, on $X$.
\end{myenuma}
Further, let $t_0,\tau\in[0,T)$ with $t_0<\tau$.
Let $f \in C([t_0,\tau],X)$ and $\mr{u}\in Y$, and define
\begin{equation}\label{C80}
  u(t) \coloneqq S(t-t_0)\mr{u} + \int\limits_{t_0}^t S(t-s)f(s)\,\rd s,
  \qquad t\in[t_0,\tau].
\end{equation}
Then $u$ satisfies
\begin{equation}\label{C81}
  u(t) = S_{\gamma}(t-t_0)\mr{u}
  + \int\limits_{t_0}^t S_{\gamma}(t-s)\bigl[f(s)+\gamma Hu(s)\bigr]\,\rd s,
  \qquad t\in[t_0,\tau].
\end{equation}
\end{lemma}

\begin{proof}
We first consider the case where $f \in C([t_0,\tau],Y)$.
From \eqref{C80} and Lemma~\ref{C77} we obtain that,
for $t\in[t_0,\tau]$,
\begin{align*}
  u(t) &= S_{\gamma}(t-t_0)\mr{u} + \int\limits_0^{t-t_0} S_{\gamma}(t-t_0-s)\gamma HS(s)\mr{u}\,\rd s
  \\
  &\quad + \int\limits_{t_0}^t \Biggl[S_{\gamma}(t-s)f(s)
  + \int\limits_0^{t-s} S_{\gamma}(t-s-r)\gamma HS(r)f(s)\,\rd r\Biggr]\rd s
  \displaybreak[0]\\[1ex]
  &= S_{\gamma}(t-t_0)\mr{u} + \int\limits_{t_0}^t S_{\gamma}(t-s)f(s)\,\rd s
  \\
  &\quad + \int\limits_0^{t-t_0} S_{\gamma}(t-t_0-s)\gamma HS(s)\mr{u}\,\rd s
  + \int\limits_{t_0}^t \int\limits_0^{t-s} S_{\gamma}(t-s-r)\gamma HS(r)f(s)\,\rd r\,\rd s.
\end{align*}
We rewrite the sum of the integrals in the last line, which we denote by $V$, as follows:
\begin{align*}
  V &= \int\limits_{t_0}^t S_{\gamma}(t-v)\gamma HS(v-t_0)\mr{u}\,\rd v
  + \int\limits_{t_0}^t \int\limits_s^t S_{\gamma}(t-v)\gamma HS(v-s)f(s)\,\rd v\,\rd s
  \\[1ex]
  &= \int\limits_{t_0}^t S_{\gamma}(t-v)\gamma HS(v-t_0)\mr{u}\,\rd v
  + \int\limits_{t_0}^t S_{\gamma}(t-v)\gamma H
  \Biggl(\,\int\limits_{t_0}^v S(v-s)f(s)\,\rd s\Biggr)\rd v
  \displaybreak[0]\\[1ex]
  &= \int\limits_{t_0}^t S_{\gamma}(t-v)\gamma H\Biggl[S(v-t_0)\mr{u}
  + \int\limits_{t_0}^v S(v-s)f(s)\,\rd s\Biggr]\rd v
  \\[1ex]
  &= \int\limits_{t_0}^t S_{\gamma}(t-v)\gamma Hu(v)\,\rd v,
\end{align*}
which shows \eqref{C81} for $f\in C([t_0,\tau],Y)$.

Now let $f\in C([t_0,\tau],X)$.  Since $C([t_0,\tau],Y)$ is dense in $C([t_0,\tau],X)$,
there exist $f_n\in C([t_0,\tau],Y)$, $n\in\NN$, such that $f_n\to f$ in $C([t_0,\tau],X)$
as $n\to\infty$.  Define
\[
  u_n(t) \coloneqq S(t-t_0)\mr{u} + \int\limits_{t_0}^t S(t-s)f_n(s)\,\rd s,
  \qquad t\in[t_0,\tau].
\]
It follows from \eqref{C59} that $u_n\to u$ in $C([t_0,\tau],Y)$ as $n\to\infty$.
Since $H\in\cB(Y,X)$, this yields that $f_n+\gamma Hu_n \to f+\gamma Hu$
in $C([t_0,\tau],X)$.  From the first part of the proof we know that
\[
  u_n(t) = S_{\gamma}(t-t_0)\mr{u} + \int\limits_{t_0}^t S_{\gamma}(t-s)
  \bigl[f_n(s)+\gamma Hu_n(s)\bigr]\,\rd s.
\]
Now we can take the limit in $X$ as $n\to\infty$ to
obtain \eqref{C81}.
\end{proof}

In the next lemma we show that the assumptions of Theorem~\ref{C31} are satisfied
if $G$ is replaced by $G-\gamma H$, so that we can apply, in particular,
Lemma~\ref{C72}\,(iii) to $G-\gamma H$.

\begin{lemma}\label{C82}
Let the assumptions of Theorem~\ref{C32} hold.
Then, for each $\gamma\ge0$, assumptions \textup{(a)} and \textup{(b)} in Theorem~\ref{C31}
are satisfied with $G$ and $(S(t))_{t\ge0}$ replaced by $G-\gamma H$ and $(S_\gamma(t))_{t\ge0}$ respectively,
and the relation
\begin{equation}\label{C83}
  \big\|S_\gamma(t)\restrict_Y\big\|_{\cB(Y)} \le \big\|S(t)\restrict_Y\big\|_{\cB(Y)},
  \qquad t\in[0,\infty),
\end{equation}
holds.
\end{lemma}

\begin{proof}
First note that $H\in\cB(Y,X)$ by \cite[Theorem~2.65]{banasiak_arlotti2006perturbations}.
Since $\gamma\ge0$, $H\ge0$, $S(s)\ge0$ and $S_\gamma(s)\ge0$ for all $s\in[0,\infty)$,
it follows from Lemma~\ref{C77} that
\[
  S(t)y \ge S_\gamma(t)y \qquad\text{for all} \;\; y\in Y_+,\; t\in[0,\infty).
\]
Hence $S(t)|_Y\ge S_\gamma(t)|_Y$ where $S(t)|_Y$ and $S_\gamma(t)|_Y$
are considered as operators from $Y$ to $Y$.  Together, with Lemma~\ref{C48},
this implies \eqref{C83}.

Relation \eqref{C83}, in combination with Lemma~\ref{C63}\,(ii), shows that
\[
  \bigl\{\big\|S_\gamma(t)\restrict_Y\big\|_{\cB(Y)}: t \in [0,t_1]\bigr\}
\]
is bounded for every $t_1\in(0,\infty)$.  It follows from Lemma~\ref{C63}\,(iii),
applied to $(S_\gamma(t))_{t\ge0}$ instead of $(S(t))_{t\ge0}$, that the function
\[
  t \mapsto \int\limits_0^t S_\gamma(t-s)\varphi(s)\,\rd s
\]
is continuous from $[0,s_1]$ to $Y$ for every $s_1>0$ and $\varphi\in C([0,s_1],X)$.
Hence we can deduce from \eqref{C78} that $t\mapsto S_\gamma(t)y$ is continuous
from $[0,\infty)$ to $Y$ for every $y\in Y$.
This shows that assumption \textup{(a)} in Theorem~\ref{C31} is satisfied
for $(S_\gamma(t))_{t\ge0}$.

Now let $s_0,s_1\in[0,\infty)$ with $s_0<s_1$, let $f\in C([s_0,s_1],Y)_+$, and set
\[
  u(t) \coloneqq \int\limits_{s_0}^t S(t-s)f(s)\,\rd s,
  \qquad t\in[s_0,s_1].
\]
Since $S(t-s)$ is positive and leaves $Y$ invariant, we obtain that $u(t)\in Y_+$
for $t\in[s_0,s_1]$.  It follows from Lemma~\ref{C79}
with $\mr u=0$ that
\[
  \int\limits_{s_0}^{s_1}S(s_1-s)f(s)\,\rd s = u(s_1)
  = \int\limits_{s_0}^{s_1}S_\gamma(s_1-s)f(s)\,\rd s
  + \int\limits_{s_0}^{s_1}S_\gamma(s_1-s)\gamma Hu(s)\,\rd s.
\]
The first integral on the right-hand side is in $Y_+$ since $S_\gamma(s_1-s)$ is
positive and leaves $Y$ invariant by assumption.
The second integral on the right-hand side is positive since $u(s)\ge0$, $H\ge0$,
$\gamma\ge0$ and $S_\gamma(s_1-s)\ge0$.  This shows that
\[
  \int\limits_{s_0}^{s_1}S(s_1-s)f(s)\,\rd s 
  \ge \int\limits_{s_0}^{s_1}S_\gamma(s_1-s)f(s)\,\rd s 
\]
for every $f\in C([s_0,s_1],Y)_+$.  Since $f\mapsto\int_{s_0}^{s_1}S(s_1-s)f(s)\rd s$
is a bounded positive operator from $C([s_0,s_1],X)$ to $Y$ by \eqref{C59},
we can apply Lemma~\ref{C48}
to deduce that the mapping $f\mapsto\int_{s_0}^{s_1}S_\gamma(s_1-s)f(s)\rd s$
can be extended to a bounded operator from $C([s_0,s_1],X)$ to $Y$ with
\[
  \Bigg\|\int\limits_{s_0}^{s_1} S_\gamma(s_1-s)\varphi(s)\,\rd s \Bigg\|_Y
  \le \Bigg\|\int\limits_{s_0}^{s_1} S(s_1-s)\varphi(s)\,\rd s \Bigg\|_Y
  \le \eta(s_1-s_0)\big\|\varphi\big\|_{C([s_0,s_1],X)}\,,
\]
which shows that assumption (b) of Theorem~\ref{C31}
is satisfied.
\end{proof}

We are now ready to prove Theorem~\ref{C32}.

\begin{proof}[Proof of Theorem~\ref{C32}]
Note again that $H\in\cB(Y,X)$ by \cite[Theorem~2.65]{banasiak_arlotti2006perturbations}.
Let $\gamma\ge0$.  Theorem~\ref{C31}\,(i)
implies that $u$ is a solution of \eqref{C57} with $t_0=0$.
Hence, by Lemma~\ref{C79} with $f(t)=F(t,u(t))$,
the function $u$ also satisfies
\[
  u(t) = S_\gamma(t)\mr u
  + \int\limits_0^t S_\gamma(t-s)\bigl[F\bigl(s,u(s)\bigr)+\gamma Hu(s)\bigr]\,\rd s,
  \qquad t\in[0,\tmax),
\]
i.e.\ it is a mild solution of the ACP \eqref{C54} and \eqref{C55}
with $t_0=0$, $G$ replaced by $G-\gamma H$, and $F$ replaced by $F_\gamma$ where
\[
  F_\gamma(t,v) \coloneqq F(t,v) + \gamma Hv,
  \qquad t\in[0,T),\; v\in Y.
\]
By assumption and Lemma~\ref{C82}, $G-\gamma H$ satisfies (a) and (b)
of Theorem~\ref{C31}.  Moreover, $F_\gamma$ satisfies condition (c) by \cref{C146}. 

Now let $\tau\in(0,\tmax)$ and $C>0$ such that $\|u(t)\|_Y\le C$ for all $t\in[0,\tau]$.
Define $\hatM^{(\gamma)}(t)$ and $\hatM_0^{(\gamma)}$ as in \eqref{C65}
and \eqref{C66} respectively, with $S(t)$ replaced by $S_\gamma(t)$.
It follows from \eqref{C83} that $\hatM^{(\gamma)}(t)\le\hatM(t)$
for $t\ge0$ and hence $\hatM_0^{(\gamma)}\le\hatM_0$.
Now choose $r>2C\hatM_0$, which then also satisfies $r>2C\hatM_0^{(\gamma)}$.
By assumption, there exists $\gamma\ge0$ such that \eqref{C62}
holds.  Now Lemma~\ref{C72}\,(iii)
implies that $u(t)\ge0$ for $t\in[0,\delta]$ with $\delta=\delta^{(\gamma)}(C,r,\tau)$,
where $\delta^{(\gamma)}(C,r,\tau)$ is defined as in \eqref{C68}
with $S(t)$ replaced by $S_\gamma(t)$.  Since $\|u(t)\|_Y\le C$ for $t\in[0,\tau]$,
we can apply that lemma repeatedly a finite number of times with the same $\delta$
(in the last step $\delta$ may be smaller) to show that $u(t)\ge0$ for $t\in[0,\tau]$.
Since $\tau\in(0,\tmax)$ was arbitrary, we obtain $u(t)\ge0$ for all $t\in[0,\tmax)$.
\end{proof}

\subsection{Special situations and additions}
\label{C15}

We start with a lemma that is useful to check assumption (e) of
Theorem~\ref{C32}; in particular, it is used in the proof of Proposition~\ref{C86} below.

\begin{lemma}\label{C87}
Let $G_0$ and $G_1$ be generators of substochastic semigroups
such that $\cD(G_0)\subseteq\cD(G_1)$ and that $G_0+G_1$ generates
a C$_0$-semigroup $(S(t))_{t\ge0}$.
Then $(S(t))_{t\ge0}$ is also substochastic.
\end{lemma}

\begin{proof}
Let $(S_i(t))_{t\ge0}$ be the semigroup generated by $G_i$, $i=0,1$.
It follows from the assumption that
\[
  \bigg\|\biggl(S_0\Bigl(\frac{t}{n}\Bigr)S_1\Bigl(\frac{t}{n}\Bigr)\biggr)^n\bigg\|
  \le 1
\]
for all $t\in[0,\infty)$ and $n\in\NN$.
Since $G_0+G_1$ generates a C$_0$-semigroup, the operator $\lambda-G_0-G_1$
is surjective for large enough $\lambda$.
Hence we can apply \cite[Corollary~III.5.8]{engel_nagel2000}, which yields that
the semigroup $(S(t))_{t\ge0}$ is given by the Trotter product formula:
\[
  S(t)x = \lim_{n\to\infty}\biggl(S_0\Bigl(\frac{t}{n}\Bigr)
  S_1\Bigl(\frac{t}{n}\Bigr)\biggr)^n x,
  \qquad t\in[0,\infty),\, x\in X.
\]
From this we easily deduce that $(S(t))_{t\ge0}$ is substochastic.
\end{proof}

\begin{proposition}\label{C86}
Let $G$ be the generator of an analytic semigroup, $(S(t))_{t\ge0}$, on a Banach space $X$.
Further, let $\alpha\in(0,1)$ and let either $Y=D_G(\alpha,p)$ with $1\le p<\infty$
or $Y=D_G(\alpha)$ or $Y=\cD((-G)^\alpha)$; 
in the case $Y=\cD((-G)^\alpha)$ assume in addition that $(S(t))_{t\ge0}$ has a strictly negative growth bound.
\begin{myenum}
\item
Then assumptions \textup{(a)} and \textup{(b)}
in Theorem~\ref{C31} are satisfied.
\item
Assume, in addition, that $X$ and $Y$ are Banach lattices such that $Y$ 
is a sublattice of $X$, assume that $G$ generates a substochastic \textup{(}analytic\textup{)} semigroup
and let $H:Y\to X$ be a positive linear operator 
so that $-H$ generates a substochastic semigroup on $X$.
Then assumptions \textup{(d)} and \textup{(e)} in Theorem~\ref{C32}
are satisfied.
\end{myenum}
\end{proposition}

\begin{proof}
(i)
It follows from \cite[Proposition~2.2.7]{lunardi1995analytic} that assumption (a)
of Theorem~\ref{C31} is satisfied.

Let us now show that assumption (b) holds.
Fix $\tau\in[0,\infty)$, and let $s_0,s_1\in[0,\tau]$ with $s_0<s_1$.
Let $\varphi\in C([s_0,s_1],X)$ and consider
\[
  v(t) = \int\limits_{s_0}^{t+s_0} S(t+s_0-s)\varphi(s)\,\rd s
  = \int\limits_0^{t}S(t-r)\varphi(r+s_0)\,\rd r,
  \qquad t\in[0,s_1-s_0].
\]
It follows from \cite[Proposition~4.2.1]{lunardi1995analytic}
(see also \cite[Lemma~7.1.1]{lunardi1995analytic})
that $v\in C^{1-\alpha}([0,s_1-s_0],Y)$ and
\[
  \|v\|_{C^{1-\alpha}([0,s_1-s_0],Y)} \le C\|\varphi\|_{C([s_0,s_1],X)}
\]
with some constant $C>0$, which depends only on $\tau$, $\alpha$ 
and the semigroup $(S(t))_{t\ge0}$.
Since $v(0)=0$, this implies that
\begin{align*}
  \Bigg\|\int\limits_{s_0}^{s_1}S(s_1-s)\varphi(s)\,\rd s\Bigg\|_Y
  &= \big\|v(s_1-s_0)\big\|_Y
  \le [v]_{C^{1-\alpha}([0,s_1-s_0],Y)}(s_1-s_0)^{1-\alpha}
  \\[1ex]
  &\le C(s_1-s_0)^{1-\alpha}\|\varphi\|_{C([s_0,s_1],X)},
\end{align*}
where we used the notation from \eqref{C95}.
With a sequence of different $\tau$ tending to $\infty$, one can construct
a function $\eta$ with the desired properties.

(ii)
It follows from \cite[Corollary~2.2.3\,(iv)]{lunardi1995analytic}
that $\cD(G)$ is dense in $Y$, and hence assumption (d) in Theorem~\ref{C32} is satisfied.

Let us now prove that assumption (e) is satisfied.
First note that $H\in\cB(Y,X)$ by \cite[Theorem~2.65]{banasiak_arlotti2006perturbations}.
It follows from \cite[Proposition~2.4.1\,(i)]{lunardi1995analytic} and Lemma~\ref{C87} that
the operator $G-\gamma H$ generates an analytic, substochastic semigroup on $X$.
Since the interpolation spaces for $G$ and $G-\gamma H$ coincide 
by \cite[Corollary~2.2.3\,(i)]{lunardi1995analytic}, the semigroup generated by $G-\gamma H$ 
leaves $Y$ invariant.
\end{proof}

\noindent
In the following proposition we prove an invariance result.  In the next section it is used 
to prove mass conservation of a solution of the coagulation--fragmentation equation.
If $X$ is an AL-space and $\phi$ is the linear extension of the norm from the positive cone,
this proposition can be used to show conservation of the norm for positive solutions.

\begin{proposition}\label{C84}
Let $X$ and $Y$ be Banach spaces such that $Y$ is continuously embedded in $X$,
let $G$ be the generator of a $C_0$-semigroup, $(S(t))_{t\ge0}$, let $T\in(0,\infty]$,
and let $F:[0,T)\times Y\to X$ be continuous.
Further, let $\phi$ be a bounded linear functional on $X$ and assume that the following
conditions are satisfied
\begin{alignat}{2}
	& \phi\bigl(S(t)f\bigr) = \phi(f) \qquad && \text{for all} \ f\in X,\, t\in[0,\infty),
	\label{C123}
	\\[1ex]
	& \phi\bigl(F(t,g)\bigr) = 0 \qquad && \text{for all} \ g\in Y,\, t\in[0,T).
	\label{C124}
\end{alignat}
Let $\mr{u}\in Y$ and let $u\in C([0,t_1),Y)$ be a mild solution of \eqref{C54}, \eqref{C55}
\textup{(}with $t_0=0$\textup{)}.  Then
\begin{equation}\label{C85}
	\phi\bigl(u(t)\bigr) = \phi(\mr{u}) \qquad \text{for all} \ t\in[0,t_1).
\end{equation}
\end{proposition}

\begin{proof}
We obtain from \eqref{C57}, \eqref{C123} and \eqref{C124} that, for $t\in[0,t_1)$,
\begin{align*}
	\phi\bigl(u(t)\bigr)
	&= \phi\bigl(S(t)\mr{u}\bigr)
	+ \phi\Biggl(\int\limits_0^t S(t-s)F\bigl(s,u(s)\bigr)\,\rd s\Biggr)
	\\[1ex]
	&= \phi(\mr{u})
	+ \int\limits_0^t \phi\Bigl(S(t-s)F\bigl(s,u(s)\bigr)\Bigr)\,\rd s
	\\[1ex]
	&= \phi(\mr{u})
	+ \int\limits_0^t \phi\Bigl(F\bigl(s,u(s)\bigr)\Bigr)\,\rd s
	= \phi(\mr{u}),
\end{align*}
which proves \eqref{C85}.
\end{proof}

\noindent
The following proposition, which is based on Gr\"onwall's inequality can be used to obtain
a priori estimates, which are used for proving global existence of solutions.

\begin{proposition}\label{C143}
Assume that \textup{(a)--(f)} in \cref{C31,C32} hold with $X=Y$ and $H=I$.
Further, let $\phi$ be a positive, bounded, linear functional on $X$ such that
\begin{equation}\label{C176}
	\phi\bigl(S(t)f\bigr) \le \phi(f) \qquad\text{for all} \ f\in X_+,\; t\ge0.
\end{equation}
Moreover, let $\mr u\in X_+$, and let $u$ be the unique mild solution on $[0,\tmax)$
of \eqref{C54}, \eqref{C55}.  Let $\tau\in(0,\tmax)$ and $\nu(\tau)\in\RR$ be such that
\begin{equation}\label{C177}
	\phi\bigl(F(t,u(t)\bigr) \le \nu(\tau)\phi\bigl(u(t)\bigr)
	\qquad\text{for all} \ t\in[0,\tau].
\end{equation}
Then
\begin{equation}\label{C175}
	\phi\bigl(u(t)\bigr) \le \phi(\mr u)e^{\nu(\tau)t}, \qquad t\in[0,\tau].
\end{equation}
\end{proposition}

\begin{proof}
Let $r\coloneq\max_{t\in[0,\tau]}\|u(t)\|_X$.  By assumption (f) there exists $\gamma\ge0$
such that \eqref{C62} holds and hence
\begin{equation}\label{C178}
	F\bigl(t,u(t)\bigr)+\gamma u(t) \ge 0 \qquad\text{for all} \ t\in[0,\tau].
\end{equation}
Since $S_\gamma(t)=e^{-\gamma t}S(t)$ with the notation from \cref{C32}, it follows from \cref{C79} that $u$ satisfies
\[
	u(t) = e^{-\gamma t}S(t)\mr u + \int\limits_0^t e^{-\gamma(t-s)}S(t-s)\bigl[F\bigl(s,u(s)\bigr)+\gamma u(s)\bigr]\,\rd s.
\]
Applying $\phi$ on both sides and using \eqref{C178}, \eqref{C176} and \eqref{C177} we obtain
\begin{align*}
	\phi\bigl(u(t)\bigr) &= e^{-\gamma t}\phi\bigl(S(t)\mr u\bigr)
	+ e^{-\gamma t}\int\limits_0^t e^{\gamma s}\phi\Bigl(S(t-s)\bigl[F\bigl(s,u(s)\bigr)+\gamma u(s)\bigr]\Bigr)\,\rd s
	\\[1ex]
	&\le e^{-\gamma t}\phi(\mr u) + e^{-\gamma t}\int\limits_0^t e^{\gamma s}\phi\Bigl(F\bigl(s,u(s)\bigr)+\gamma u(s)\Bigr)\,\rd s
	\\[1ex]
	&\le e^{-\gamma t}\phi(\mr u) + e^{-\gamma t}\bigl(\nu(\tau)+\gamma\bigr)\int\limits_0^t e^{\gamma s}\phi\bigl(u(s)\bigr)\,\rd s.
\end{align*}
Setting $y(t)=e^{\gamma t}\phi(u(t))$ we arrive at
\[
	y(t) \le \phi(\mr u) + \bigl(\nu(\tau)+\gamma\bigr)\int\limits_0^t y(s)\,\rd s,
\]
and Gr\"onwall's inequality yields
\[
	y(t) \le \phi(\mr u)e^{(\nu(\tau)+\gamma)t},
\]
which, in turn, implies \eqref{C175}.
\end{proof}

\section{The Coagulation--Fragmentation System}
\label{C04}

In this section we apply the abstract results from Section~\ref{C03} to the C--F system, 
posed as a semi-linear ACP in an $\ell_w^1$ space, and prove our main result, \cref{C100}.
We begin this section by imposing assumptions on the weights $(w_n)_{n=1}^{\infty}$, the fragmentation coefficients $b_{n,j}$
and on the coagulation rates $k_{n,j}$ (see \cref{C88,C98}), and obtain numerous useful properties under these assumptions, 
including the Lipschitz continuity of an operator used to describe the coagulation terms in the C--F system. 
We then pose the C--F system as a semi-linear ACP and apply the results from Section~\ref{C03} 
to prove existence, uniqueness and positivity of solutions. 
In particular, we apply these abstract results in two settings, corresponding to different choices of the space $Y$. 
The first setting corresponds to $Y=X$, while the second setting considers the case 
where $Y$ is an interpolation space between $X$ and an appropriate dense subspace of $X$. 
The second setting allows solutions to be obtained under weaker assumptions on the coagulation rates than are required in the $Y=X$ setting.
Further, in \cref{C158} we show that, under extra assumptions, the mild solutions are classical solutions.
Finally, in \cref{C152} we establish global existence of solutions in a particular situation.

\subsection{Setting of the Problem}\label{C150}

In the following two subsections (\S\ref{C153}, \ref{C154}) we consider operators corresponding to the fragmentation terms
and to the coagulation terms separately.  
In \cref{C153} we recall some facts from \cite{kerr_lamb_langer2020fragpaper} about the generator of a semigroup
that corresponds to the fragmentation terms.  In \cref{C154} we then study the non-linear coagulation operator.
In \cref{C155} we provide some sufficient conditions for \cref{C88,C98}.

\subsubsection{The Fragmentation Operator}\label{C153}

Recall that, throughout the paper, we let Assumption~\ref{C22} hold.  Moreover, we also make the 
following assumption on the weight $w=(w_n)_{n=1}^{\infty}$ and its relation to the 
fragmentation coefficients $b_{n,j}$.

\begin{assumption}\label{C88}
Let $w=(w_n)_{n=1}^\infty$ be such that $w_n \ge n$ for all $n \in \mathbb{N}$, $(w_n)_{n=1}^{\infty}$ 
is monotone increasing and there exists $\kappa \in (0,1]$ such that
\begin{equation}\label{C126}
	\sum\limits_{n=1}^{j-1} w_nb_{n,j} \le \kappa w_j \qquad \text{for all} \ j\in\{2,3,\ldots\}.
\end{equation}
\hfill $\lozenge$
\end{assumption}

\medskip

\noindent
We provide sufficient conditions for \cref{C88} in \cref{C186}.
Note that, under \cref{C88} we can conclude that $\ell_w^1$ is continuously embedded in $X_{[1]}$, 
where $X_{[1]}$ is defined in \eqref{C147}.

In \cref{C98} and \cref{C100} it is important to distinguish the two cases $\kappa=1$
and $\kappa<1$.  
The latter case, which is more restrictive on the fragmentation coefficients $b_{n,j}$,
allows faster growing coagulation coefficients.

\medskip

Let us introduce the following formal expressions, which correspond to the two
linear fragmentation terms on the right-hand side of \eqref{C21},
\begin{align*}
	\mathcal{A}: (f_n)_{n=1}^{\infty} \mapsto \bigl(-a_nf_n\bigr)_{n=1}^{\infty}
	\qquad \text{and} \qquad
	\mathcal{B}: (f_n)_{n=1}^{\infty} \mapsto
	\Biggl(\sum\limits_{j=n+1}^{\infty} a_j b_{n,j}f_j\Biggr)_{n=1}^{\infty}.
\end{align*}
Operator realisations, $A^{(w)}$ and $B^{(w)}$, of $\mathcal{A}$ and $\mathcal{B}$
respectively, are defined in $\ell_w^1$ with their maximal domains by
\begin{alignat*}{2}
	A^{(w)}f &= \mathcal{A}f, \qquad &
	\mathcal{D}(A^{(w)}) &= \bigl\{f \in \ell_w^1: \mathcal{A}f \in \ell_w^1\bigr\},
	\\[1ex]
	B^{(w)}f &= \mathcal{B}f, \qquad &
	\mathcal{D}(B^{(w)}) &= \bigl\{f \in \ell_w^1: \mathcal{B}f \in \ell_w^1\bigr\}.
\end{alignat*}
It is proved in \cite[Theorem~3.4]{kerr_lamb_langer2020fragpaper} that
\begin{equation}\label{C127}
	G^{(w)} \coloneq \overline{A^{(w)}+B^{(w)}},
\end{equation}
where the overline denotes the closure,
is the generator of a substochastic semigroup, $(S^{(w)}(t))_{t\ge0}$, on $\ell_w^1$.
If $\kappa<1$ in \eqref{C126}, then the closure in \eqref{C127} is not needed, and
$G^{(w)}=A^{(w)}+B^{(w)}$ generates an analytic semigroup; see \cite[Theorem~5.2]{kerr_lamb_langer2020fragpaper}.

\subsubsection{The Coagulation Operator}\label{C154}

Let us now state the assumptions on the coagulation coefficients $k_{n,j}$ that are used in \cref{C100};
we consider two situations; see (CI) and (CII) below.

\begin{assumption}\label{C98}
Let \cref{C22,C88} be satisfied, let $T\in(0,\infty]$ and $\alpha\in[0,1)$.
Assume that one of the following conditions is satisfied:
\begin{enumerate}[label=\textup{(C\Roman*)}]
\item\label{C128}
	$\kappa\in(0,1]$ and $\alpha=0$;
\item\label{C129}
	$\kappa\in(0,1)$ and $\alpha\in(0,1)$.
\end{enumerate}
Further, assume that
\begin{enumerate}[label=(A\arabic*)]
\item
	$k_{n,j}:[0,T)\to[0,\infty)$ are continuous for all $n,j\in\NN$;
\item
	for every $t'\in(0,T)$ there exists a constant $c(t')>0$ such that
	\begin{equation}\label{C116}
		k_{n,j}(t) \le c(t')\frac{\tildew_n\tildew_j}{w_{n+j}}
		\qquad\text{for all} \ n,j\in\NN \ \text{and} \ t\in[0,t'],
	\end{equation}
	where
	\begin{equation}\label{C125}
		\tildew_n \coloneq (1+a_n)^\alpha w_n, \qquad n\in\NN.
	\end{equation}
	Further, we set $\tildew=(\tildew_n)_{n=1}^\infty$.
\hfill $\lozenge$
\end{enumerate}
\end{assumption}

\medskip

\noindent
In \cref{C187} (in \S\ref{C155}) we present some sufficient conditions for \cref{C98}.

In the following proposition we define a non-linear operator that corresponds to the coagulation terms
in the second line of \eqref{C21} and prove some properties.
This proposition is an extension of \cite[Theorem~4.3\,(ii)]{mcbride_smith_lamb2010strongly}.

\begin{proposition}\label{C115}
Let \cref{C22,C98} hold and let $\tildew$ be as in \eqref{C125}.
Then the operator $K^{(\tildew,w)}$ defined by
\begin{equation}\label{C131}
\begin{aligned}
	K^{(\tildew,w)}(t,f) \coloneq \biggl(\frac{1}{2}\sum_{j=1}^{n-1}k_{n-j}(t)f_{n-j}f_j 
	- \sum_{j=1}^\infty k_{n,j}(t)f_n f_j\biggr)_{n=1}^\infty, \hspace*{15ex} &
	\\[1ex]
	\qquad t\in[0,T),\; f=(f_j)_{j=1}^\infty \in\ell_{\tildew}^1, &
\end{aligned}
\end{equation}
is a continuous mapping from $[0,T) \times \ell_{\tildew}^1$ to $\ell_w^1$.  
Moreover, $K^{(\tildew,w)}$ is Lipschitz in the second argument on bounded sets, uniformly in the 
first argument on compact intervals.
\end{proposition}

\noindent
Before we prove \cref{C115} we need two lemmas.

\begin{lemma}\label{C135}
Let $\omega_n>0$ and $\beta_{n,j}\in\RR$ for $n,j\in\NN$.  Then
\begin{equation}\label{C136}
	\sum_{n=1}^\infty \omega_n\sum_{j=1}^{n-1}\beta_{n-j,j}
= \sum_{n=1}^\infty \sum_{j=1}^\infty \omega_{n+j}\beta_{n,j}
\end{equation}
provided the series on the right-hand side converges absolutely.
\end{lemma}

\begin{proof}
Interchanging the order of summation and substituting $l=n-j$ we obtain
\begin{align*}
	\sum_{n=1}^\infty \omega_n\sum_{j=1}^{n-1}\beta_{n-j,j}
	= \sum_{j=1}^\infty \sum_{n=j+1}^\infty \omega_n\beta_{n-j,j}
	= \sum_{j=1}^\infty \sum_{l=1}^\infty \omega_{l+j}\beta_{l,j},
\end{align*}
which proves \eqref{C136}.
\end{proof}

\begin{lemma}\label{C110}
Let $(w_n)_{n=1}^{\infty}$ be such that $w_n>0$ for all $n\in\NN$ 
and $(w_n)_{n=1}^{\infty}$ is monotone increasing.  
Let $\varphi_{n,j} \in \mathbb{R}$ for all $n,j\in\NN$ and set
\begin{equation}\label{C132}
	\cL_\varphi[f,g] \coloneq \biggl(\frac{1}{2}\sum\limits_{j=1}^{n-1} \varphi_{n-j,j}f_{n-j}g_j
	-\sum\limits_{j=1}^{\infty} \varphi_{n,j}f_ng_j\biggr)_{n=1}^\infty
\end{equation}
for all $f=(f_j)_{j=1}^\infty$, $g=(g_j)_{j=1}^\infty$ such that the infinite series
in \eqref{C132} converges absolutely.
Assume that
\[
	\sum\limits_{n=1}^{\infty}\sum\limits_{j=1}^{\infty} w_{n+j}|\varphi_{n,j}|\,|f_n|\,|g_j| < \infty.
\]
Then $\cL_\varphi[f,g]$ is well defined, $\cL_\varphi[f,g]\in\ell_w^1$ and
\begin{equation}\label{C112}
  \big\|\cL_\varphi[f,g]\big\|_w
  \le \frac{3}{2} \sum\limits_{n=1}^{\infty}\sum\limits_{j=1}^{\infty} w_{n+j}|\varphi_{n,j}|\,|f_n|\,|g_j|.
\end{equation}
\end{lemma}

\begin{proof}
Let us consider the two sums in \eqref{C132} separately.
For the first sum we use \cref{C135} with $\omega_n=w_n$ and $\beta_{n,j}=|\varphi_{n,j}f_n g_j|$
to obtain
\begin{equation}\label{C113}
\begin{aligned}
	\Biggl\|\biggl(\frac{1}{2}\sum\limits_{j=1}^{n-1} \varphi_{n-j,j}f_{n-j}g_j\biggr)_{n\in\NN}\Biggr\|_w
	&\le \frac{1}{2}\sum\limits_{n=1}^{\infty}w_n\sum\limits_{j=1}^{n-1} |\varphi_{n-j,j}|\,|f_{n-j}|\,|g_j|
	\\[1ex]
	&=\frac{1}{2}\sum\limits_{n=1}^{\infty}\sum\limits_{j=1}^{\infty} w_{n+j}|\varphi_{n,j}|\,|f_n|\,|g_j|.
\end{aligned}
\end{equation}
Let us now consider the second sum; the monotonicity of $(w_n)_{n=1}^{\infty}$ yields
\begin{equation}\label{C114}
\begin{aligned}
	\Bigg\|\biggl(\sum\limits_{j=1}^{\infty} \varphi_{n,j}f_n g_j\biggr)_{n\in\NN}\Biggr\|_w
	&\le \sum\limits_{n=1}^{\infty} w_n \sum\limits_{j=1}^{\infty} |\varphi_{n,j}|\,|f_n|\,|g_j|
	\\[1ex]
	&\le \sum\limits_{n=1}^{\infty} \sum\limits_{j=1}^{\infty} w_{n+j}|\varphi_{n,j}|\,|f_n|\,|g_j|.
\end{aligned}
\end{equation}
Combining \eqref{C113} and \eqref{C114} we arrive at \eqref{C112}.
\end{proof}

\noindent
In order to prove \cref{C115} we introduce a bilinear mapping that corresponds to $K^{(\tildew,w)}$:
\begin{equation}\label{C162}
	\widetilde K^{(\tildew,w)}[t,f,g] \coloneq \biggl(\frac{1}{2}\sum_{j=1}^{n-1} k_{n-j,j}(t)f_{n-j}g_j
	- \sum_{j=1}^{\infty} k_{n,j}(t)f_n g_j\biggr)_{n=1}^\infty,
\end{equation}
which, as we see below, is well defined for $f,g\in\ell_{\tildew}^1$.

\begin{proof}[Proof of \cref{C115}]
Let $t'\in(0,T)$.  Then there exists $c(t')>0$ such that \eqref{C116} holds.  With $\varphi_{n,j}=k_{n,j}(t)$ we
obtain from \cref{C110} that, for $f,g\in\ell_{\tildew}^1$ and $t\in[0,t']$,
\begin{align*}
  \big\|\widetilde K^{(\tildew,w)}[t,f,g]\big\|_w &= \big\|\cL_\varphi[f,g]\big\|_w
  \le \frac{3}{2}\sum_{n=1}^\infty\sum_{j=1}^\infty w_{n+j}|k_{n,j}(t)|\,|f_n|\,|g_j|
  \\[1ex]
  &\le \frac{3}{2}\sum_{n=1}^\infty\sum_{j=1}^\infty c(t')\tildew_n\tildew_j|f_n|\,|g_j|
  = \frac{3}{2}c(t')\|f\|_{\tildew}\|g\|_{\tildew}.
\end{align*}
This shows that $\widetilde K^{(\tildew,w)}$ is a well-defined mapping from $[0,T)\times\ell_{\tildew}^1\times\ell_{\tildew}^1$
into $\ell_w^1$ and that (a) and (b) in \cref{C29} with $\widetilde F=\widetilde K^{(\tildew,w)}$ are satisfied.
Hence $K^{(\tildew,w)}$ is Lipschitz in the second argument on bounded sets, uniformly 
in the first argument on compact intervals by \cref{C29}\,(i).

To show continuity in $t$, let $t_0\in[0,T)$ and choose $t'\in(t_0,T)$.  
For $t\in[0,t']$ and $f,g\in\ell_{\tildew}^1$ we obtain
from \cref{C110} with $\varphi_{n,j}=k_{n,j}(t)-k_{n,j}(t_0)$ that
\begin{align}
	& \big\|\widetilde K^{(\tildew,w)}[t,f,g]-\widetilde K^{(\tildew,w)}[t_0,f,g]\big\|_w
	= \big\|\cL_\varphi[f,g]\big\|_w
	\nonumber\\[1ex]
	&\le \frac{3}{2}\sum_{n=1}^\infty\sum_{j=1}^\infty w_{n+j}\big|k_{n,j}(t)-k_{n,j}(t_0)\big|\,|f_n|\,|g_j|.
	\label{C118}
\end{align}
The expression in the latter sum satisfies
\begin{equation}\label{C117}
	w_{n+j}\big|k_{n,j}(t)-k_{n,j}(t_0)\big|\,|f_n|\,|g_j| \le 2c(t')\tildew_n\tildew_j|f_n|\,|g_j|
\end{equation}
by \eqref{C116}, and the double-sum over the right-hand side of \eqref{C117} is finite.  
As $t\to t_0$, the expressions in the sum in \eqref{C118} converge to $0$ for all $n,j\in\NN$.
Therefore also the sums in \eqref{C118} converge to $0$ by the Dominated Convergence Theorem.
This shows that $t\mapsto \widetilde K^{(\tildew,w)}[t,f,g]$ is continuous from $[0,T)$ to $\ell_w^1$
for all $f,g\in\ell_{\tildew}^1$.
It follows from \cref{C29}\,(ii) that $K^{(\tildew,w)}$ is continuous from $[0,T)\times\ell_{\tildew}^1$
to $\ell_w^1$.
\end{proof}

Later we also need the following lemma.

\begin{lemma}\label{C134}
Let $\omega_n>0$ for $n\in\NN$ and $f\in\ell_{\tildew}^1$.
Further, assume that
\[
	\sum_{n=1}^\infty\sum_{j=1}^\infty \omega_{n+j}k_{n,j}(t)f_n f_j
\]
converges absolutely.
Then
\begin{equation}\label{C148}
	\sum_{n=1}^\infty \omega_n \bigl[K^{(\tildew,w)}(t,f)\bigr]_n
	= \frac{1}{2}\sum_{n=1}^\infty\sum_{j=1}^\infty\bigl(\omega_{n+j}-\omega_n-\omega_j\bigr)
	k_{n,j}(t)f_n f_j.
\end{equation}
\end{lemma}

\begin{proof}
From \cref{C135} with $\beta_{n,j}=k_{n,j}(t)f_n f_j$ we obtain
\begin{align*}
	& \sum_{n=1}^\infty \omega_n \bigl[K^{(\tildew,w)}(t,f)\bigr]_n
	= \frac{1}{2}\sum_{n=1}^\infty \omega_n\sum_{j=1}^{n-1} k_{n-j,j}(t)f_{n-j}f_j
	- \sum_{n=1}^\infty \omega_n\sum_{j=1}^\infty k_{n,j}(t)f_n f_j
	\\[1ex]
	&= \frac{1}{2}\sum_{n=1}^\infty\sum_{j=1}^\infty \omega_{n+j}k_{n,j}(t)f_n f_j
	- \sum_{n=1}^\infty\sum_{j=1}^\infty \omega_n k_{n,j}(t)f_n f_j
	\\[1ex]
	&= \frac{1}{2}\sum_{n=1}^\infty\sum_{j=1}^\infty \omega_{n+j}k_{n,j}(t)f_n f_j
	- \frac{1}{2}\sum_{n=1}^\infty\sum_{j=1}^\infty \omega_n k_{n,j}(t)f_n f_j
	- \frac{1}{2}\sum_{j=1}^\infty\sum_{n=1}^\infty \omega_n k_{j,n}(t)f_n f_j,
\end{align*}
which proves \eqref{C148}.
\end{proof}

\subsubsection{Sufficient Conditions}\label{C155}

In this section we consider some sufficient conditions for \cref{C88,C98}.
Let us start with \cref{C88}.

\pagebreak[3]

\begin{proposition}\label{C186}
\rule{0ex}{1ex}
\begin{myenum}
\item
	If \eqref{C149} is satisfied, then \eqref{C126} holds with $\kappa=1$ and $w_n=n^p$ 
	for all $p\ge1$.
\item
	For any $b_{n,j}\ge0$, one can always find a weight $w$ such that \eqref{C126} holds with $\kappa<1$.
\item
	If \eqref{C149} is satisfied, then one can choose $w_n=r^n$ with $r>2$
	such that \eqref{C126} holds with $\kappa=\frac{2}{r}<1$.
\item
	Let $\phi:(0,1)\to[0,\infty)$ be a function such that $x\mapsto x\phi(x)$ is
	non-decreasing and integrable on $(0,1)$.
	Define
	\[
		b_{n,j} \coloneq d_j \phi\Bigl(\frac{n}{j}\Bigr), \qquad n,j\in\NN,\,n<j,
	\]
	where
	\[
		d_j = \frac{j}{\sum\limits_{m=1}^{j-1}m\phi\bigl(\frac{m}{j}\bigr)}.
	\]
	Then there exist $p>1$ and $\kappa<1$ such that \eqref{C126} holds with $w_n=n^p$.
\end{myenum}
\end{proposition}

\begin{proof}
The statement in (i) follows from \cite[Remark~3.3]{kerr_lamb_langer2020fragpaper} 
since $\bigl(\frac{n^p}{n}\bigr)_{n=1}^{\infty}$ is monotone increasing.
Items (ii) and (iii) follow directly from \cite[Theorem~5.5]{kerr_lamb_langer2020fragpaper}.

Finally, we come to the proof of (iv).  Let $p>1$.  Then also $x\mapsto x^p\phi(x)$ is non-decreasing
and integrable on $(0,1)$.  We have to estimate the following expression from above:
\[
	\frac{1}{w_j}\sum_{n=1}^{j-1}w_n b_{n,j} = \frac{1}{j^p}\sum_{n=1}^{j-1} d_j n^p\phi\Bigl(\frac{n}{j}\Bigr)
	= \frac{1}{j^{p-1}}
	\cdot\frac{\sum\limits_{n=1}^{j-1}n^p\phi\bigl(\frac{n}{j}\bigr)}{\sum\limits_{n=1}^{j-1}n\phi\bigl(\frac{n}{j}\bigr)}.
\]
Let us deal with numerator and denominator separately: for $j\ge2$ we have
\begin{align*}
	\sum_{n=1}^{j-1}n^p\phi\Bigl(\frac{n}{j}\Bigr)
	&\le \int_1^j x^p\phi\Bigl(\frac{x}{j}\Bigr)\,\rd x
	\le \int_0^j x^p\phi\Bigl(\frac{x}{j}\Bigr)\,\rd x
	= j^{p+1}\int_0^1 t^p\phi(t)\,\rd t,
	\\[1ex]
	\sum_{n=1}^{j-1}n\phi\Bigl(\frac{n}{j}\Bigr)
	&\ge \int_0^{j-1} x\phi\Bigl(\frac{x}{j}\Bigr)\,\rd x
	= j^2\int_0^{1-\frac{1}{j}} t\phi(t)\,\rd t
	\ge j^2\int_0^{\frac12} t\phi(t)\,\rd t,
\end{align*}
which then yields
\begin{equation}\label{C188}
	\frac{1}{w_j}\sum_{n=1}^{j-1}w_n b_{n,j}
	\le \frac{\int_0^1 t^p\phi(t)\,\rd t}{\int_0^{\frac12} t\phi(t)\,\rd t}.
\end{equation}
The Monotone Convergence Theorem implies that the numerator converges to $0$ as $p\to\infty$.
Hence there exists $p>1$ such that the right-hand side of \eqref{C188} is strictly less than $1$.
With this choice of $p$ the inequality in \eqref{C126} is satisfied with $\kappa$ equal to 
the right-hand side of \eqref{C188}.
\end{proof}

\begin{remark}\label{C182}
\rule{0ex}{1ex}
\begin{myenum}
\item
Note that the coefficients $b_{n,j}$, as defined in \cref{C186}\,(iv), satisfy the second relation in \eqref{C26}.
Moreover, one can choose, e.g.\ $\phi(x)=x^\nu$ with $\nu\ge-1$, which leads to the special case
that is considered in \cite[Example~2.3]{kerr_lamb_langer2024nonautonomous}.
\item
The form of $b_{n,j}$ considered in \cref{C186}\,(iv) is also examined in \cite[Example~5.1.66]{banasiak_lamb_laurencot2020vol1}, 
where it is shown that, without the monotonicity assumption on $x\mapsto x\phi(x)$, the weights $w_n=n^p$, $p>1$, 
satisfy a lim\,sup condition instead of the pointwise condition \eqref{C126}. 
Moreover, for continuous fragmentation equations, \emph{homogeneous} kernels, i.e.\ kernels of the form
\[
	b(x,y) = \frac{1}{y}h\Bigl(\frac{x}{y}\Bigr), \qquad x<y,
\]
with a function $h$ that satisfies $\int_0^1 th(t)\,\rd t=1$, have also been 
considered in the literature; see, e.g.\ \cite{cheng_redner1988,edwards_cai_han1990,cai_edwards_han1991,banasiak_lamb2012analytic}.
\item
It is not always possible to use $w_n=n^p$ so that \eqref{C126} holds with $\kappa<1$.
Let
\[
	a_n = 
	\begin{cases} 
		0 & \text{if} \ n=1, \\[0.5ex] 
		n & \text{if} \ n\ge2, 
	\end{cases}
	\qquad
	b_{n,j} = 
	\begin{cases} 
		2 & \text{if} \ j=2,\,n=1, \\[0.5ex]
		1 & \text{if} \ j\ge 3,\,n\in\{1,j-1\}, \\[0.5ex]
		0 & \text{if} \ j\ge 3,\,n\in\{2,\ldots,j-2\},
	\end{cases}
\]
and $w_n=n^p$ with $p\ge1$.
In \cite[Theorem~3]{banasiak2011irregular} and \cite[Proposition~A.3]{banasiak2012global} it is shown that 
the semigroup generated by $G^{(w)}$, defined in \eqref{C127}, is not analytic in $\ell_w^1$ and 
hence \eqref{C126} cannot hold with $\kappa<1$.
However, by \cref{C186}\,(iii) we can choose $w_n=r^n$ with $r>1$ such that \eqref{C126} holds with $\kappa<1$.
Note that the coefficients $b_{n,j}$ model the situation where particles of size $j$
fragment into particles of sizes $j-1$ and $1$.
Together with an appropriate coagulation kernel, this is the model that
was studied by Becker and D\"oring in \cite{becker_doring1935}.
\hfill $\lozenge$
\end{myenum}
\end{remark}

\bigskip

\noindent
Let us now consider sufficient conditions for \cref{C98}.

\begin{proposition}\label{C187}
Let $\alpha\in[0,1)$.
\begin{myenum}
\item
Let $w_n=n^p$ with some $p \ge 1$ and assume that, for each $t'\in[0,T)$, there exists 
a constant $\hat{c}(t')>0$ such that
\begin{equation}\label{C180}
	k_{n,j}(t) \le \hat{c}(t')\bigl[\min\{n,j\}\bigr]^p(1+a_n)^{\alpha}(1+a_j)^{\alpha}
	\qquad\text{for all} \ n,j\in\NN,\;t\in[0,t'].
\end{equation}
Then \eqref{C116} is satisfied with $c(t')=2^p\hat c(t')$.
\item
Let $w_n=r^n$ with $r>1$ and assume that, for each $t'\in[0,T)$, there exists 
a constant $\hat{c}(t')>0$ such that
\[
	k_{n,j}(t) \le \hat{c}(t')(1+a_n)^{\alpha}(1+a_j)^{\alpha}
	\qquad\text{for all} \ n,j\in\NN,\;t\in[0,t'].
\]
Then \eqref{C116} is satisfied with $c(t')=\hat c(t')$.
\end{myenum}
\end{proposition}

\begin{proof}
We only prove (i); the proof for (ii) is similar but simpler.
Since
\[
	\frac{nj}{n+j} = \frac{1}{\frac{1}{j}+\frac{1}{n}}
	\ge \frac{1}{\,\frac{2}{\min\{n,j\}}\,} = \frac{\min\{n,j\}}{2},
\]
we obtain, for $t\in[0,t']$,
\begin{align*}
	k_{n,j}(t) &\le \hat{c}(t')\bigl[\min\{n,j\}\bigr]^p(1+a_n)^{\alpha}(1+a_j)^{\alpha}
	\\[1ex]
	&\le \hat{c}(t')2^p\frac{n^pj^p}{(n+j)^p}(1+a_n)^{\alpha}(1+a_j)^{\alpha},
\end{align*}
which proves \eqref{C116}.
\end{proof}

\begin{remark}\label{C185}
The theory of interpolation spaces is used by Banasiak in \cite{banasiak2019survey} to examine \eqref{C21}, 
where the fragmentation coefficients $b_{n,j}$ are such that Assumption~\ref{C88} 
holds for $w_n=n^p$ for some $p \geq 1$ and all $n \in \mathbb{N}$.  
In \cite[Theorem~4.12]{banasiak2019survey} the situation is considered when $k_{n,j}(t) \equiv k_{n,j}$  is independent 
of time and satisfies
\begin{equation}\label{C179}
	k_{n,j} \le C\bigl((1+a_n)^{\alpha}+(1+a_j)^{\alpha}\bigr),
\end{equation}
with some $C>0$ and $\alpha\in(0,1)$.
Since
\begin{align*}
	(1+a_n)^{\alpha}(1+a_j)^{\alpha}
	&= \frac{1}{2}\bigl[(1+a_n)^{\alpha}(1+a_j)^{\alpha}+(1+a_n)^{\alpha}(1+a_j)^{\alpha}\bigr]
	\\[1ex]
	&\ge \frac{1}{2}\bigl[(1+a_n)^{\alpha}+(1+a_j)^{\alpha}\bigr],
\end{align*}
we have
\begin{equation}\label{C181}
	(1+a_n)^{\alpha}+(1+a_j)^{\alpha} \le 2\bigl[\min\{n,j\}\bigr]^p(1+a_n)^{\alpha}(1+a_j)^{\alpha},
\end{equation}
for all $p \ge 1$.
Hence \eqref{C180} is much less restrictive than \eqref{C179} and, moreover, allows time-dependent coagulation coefficients.
\hfill $\lozenge$
\end{remark}

\subsection{Local Existence and Uniqueness Results of Mild Solutions}\label{C151}

Let \cref{C22,C88,C98} hold and let $G^{(w)}$, $S^{(w)}$ and $K^{(\tildew,w)}$
be as in \cref{C150}.
We prove local existence and uniqueness results for the ACP
\begin{align}
	u'(t) &= G^{(w)}u(t)+K^{(\tildew,w)}\bigl(t,u(t)\bigr), \qquad t \in (0,T), 
	\label{C137}
	\\[1ex]
	u(0) &= \mr{u}.
	\label{C138}
\end{align}
Mild solutions of \eqref{C137}, \eqref{C138} on the interval $[0,t_1)$ with $t_1\in(0,T]$ are solutions of
\begin{equation}\label{C139}
	u(t) = S^{(w)}(t)\mathring{u}+\int\limits_0^t S^{(w)}(t-s)K^{(\tildew,w)}(s,u(s))\,\rd s,
	\qquad  t \in [0,t_1);
\end{equation}
see \cref{C56}.

The following theorem, the main result of the paper, provides existence and uniqueness of solutions of \eqref{C137}, \eqref{C138}.

\begin{theorem}\label{C100}
Let \cref{C22,C88,C98} be satisfied 
and let $\mr u=(\mr u_n)_{n=1}^\infty\in(\ell_{\tildew}^1)_+$.
Then the following statements hold.
\begin{myenum}
\item
There exists $\tmax\in(0,T]$ such that \eqref{C137}, \eqref{C138} has a unique maximal positive mild solution
$u\in C([0,\tmax),\ell_{\tildew}^1)$.
\item
If $\tmax<T$, then $\|u(t)\|_{\tildew}\to\infty$ as $t\to\tmax^-$.
\item
For every $\tau\in(0,\tmax)$ there exist $r_\tau>0$ and $M_\tau>0$ such that,
whenever $\mr v=(\mr v_n)_{n=1}^\infty\in(\ell_{\tildew}^1)_+$ with $\|\mr v-\mr u\|_{\tildew}\le r_\tau$
and $v$ is the corresponding maximal mild solution,
then $v$ is defined on $[0,\tau]$ and
\[
  \|v(t)-u(t)\|_{\tildew} \le M_\tau\|\mr v-\mr u\|_{\tildew}
  \qquad\text{for all} \ t\in[0,\tau].
\]
\item
Let $\totmass$ be defined as in \eqref{C159}.
If \eqref{C149} holds, then
\begin{equation}\label{C165}
	\totmass\bigl(u(t)\bigr) \le \totmass(\mr{u}) \qquad\text{for all} \ t\in[0,\tmax).
\end{equation}
If \eqref{C26} holds, then
\begin{equation}\label{C189}
	\totmass\bigl(u(t)\bigr) = \totmass(\mr{u}) \qquad\text{for all} \ t\in[0,\tmax).
\end{equation}
\end{myenum}
\end{theorem}

\medskip

\noindent
The proof of \cref{C100} is given \cref{C157}.  Before we come to this, we need some preliminary results
in the case when \ref{C129} holds.

\subsubsection{The Setting in the Case \ref{C129}}\label{C156}

In case \ref{C129} in \cref{C98} we rewrite \eqref{C139} slightly in order to have a semigroup with 
a strictly negative growth bound.
To this end, define the operators $A_{w,1}\coloneq A^{(w)}-I$, i.e.\
\[
	A_{w,1}f = \bigl(-(a_n+1)f_n\bigr)_{n=1}^\infty, \qquad f=(f_j)_{j=1}^\infty\in\ell_w^1,
\]
and
\begin{equation}\label{C167}
	G_{w,1} \coloneq A_{w,1}+B^{(w)} = G^{(w)}-I.
\end{equation}
It is proved in \cite[Theorem~5.2]{kerr_lamb_langer2020fragpaper} that $G^{(w)}$ generates
a substochastic, analytic semigroup, and hence $G_{w,1}$ generates an analytic semigroup with growth bound at most $-1$.

We rewrite \eqref{C137} as
\[
	u'(t) = G_{w,1}u(t)+K^{(\tildew,w)}\bigl(t,u(t)\bigr)+u(t), \qquad t \in (0,T), 
\]
which corresponds to the integral equation
\begin{equation}\label{C140}
	u(t) = S_{w,1}(t)\mr{u}+\int\limits_0^t S_{w,1}(t-s)\bigl[K^{(\tildew,w)}(s,u(s))+u(s)\bigr]\,\rd s,
	\qquad  t \in [0,t_1).
\end{equation}
We show below that \eqref{C139} and \eqref{C140} are equivalent.

For $Y$ we use an interpolation space.
The following result gives a characterisation of the interpolation space $D_{G_{w,1}}(\alpha,1)$,
which is needed in the proof of \cref{C100}.

\begin{proposition}\label{C89}
Let \cref{C22,C88,C98} be satisfied.  Further, assume that \ref{C129} holds.
Then
\begin{equation}\label{C90}
	D_{G_{w,1}}(\alpha,1) = \cD\bigl((-A_{w,1})^{\alpha}\bigr) =\ell_{\tildew}^1,
\end{equation}
with equivalent norms, where $D_{G_{w,1}}(\alpha,1)$ is the interpolation space defined in \eqref{C49} 
and $\tildew$ is as in \eqref{C125}.
\end{proposition}

\begin{proof}
Let us first determine the interpolation space between $\cD(A_{w,1})$ and $\ell_w^1$, namely we prove that
\begin{equation}\label{C133}
	D_{A_{w,1}}(\alpha,1) = \ell_{\tildew}^1.
\end{equation}
By definition (see \eqref{C49}) we have
\[
	D_{A_{w,1}}(\alpha,1)
	= \Bigl\{f \in \ell_w^1: t \mapsto v(t) \coloneqq \big\|t^{-\alpha}A_{w,1}e^{A_{w,1}t}f\big\|_w \in L^1(0,1)\Bigr\}.
\]
Let $f \in \ell_{\tildew}^1$.  With the substitution $\tau=(a_n+1)t$ we obtain
\begin{align*}
	\|v(t)\|_{L^1(0,1)} &= \int\limits_0^1 \sum\limits_{n=1}^{\infty} w_nt^{-\alpha}(a_n+1)e^{-(a_n+1)t}|f_n|\,\rd t
	\\
	&= \sum\limits_{n=1}^{\infty} w_n(a_n+1)|f_n|\int\limits_0^{a_n+1}\Bigl(\frac{\tau}{a_n+1}\Bigr)^{-\alpha} e^{-\tau}\,
	\frac{\rd\tau}{a_n+1}
	\\
	&\le \sum\limits_{n=1}^{\infty} \tildew_n|f_n|\int\limits_0^{\infty}\tau^{-\alpha}e^{-\tau}\,\rd\tau
	= \Gamma(1-\alpha)\sum\limits_{n=1}^{\infty} \tildew_n|f_n| < \infty,
\end{align*}
and hence $\ell_{\tildew}^1 \subseteq D_{A_{w,1}}(\alpha,1)$.

On the other hand, \cref{C52}\,(iii) implies
\[
	D_{A_{w,1}}(\alpha,1) \subseteq \cD\bigl((-A_{w,1})^{\alpha}\bigr) = \ell_{\tildew}^1,
\]
which shows the reverse inclusion, and therefore \eqref{C133} holds.

Next, we obtain from \cite[(3.6)]{kerr_lamb_langer2020fragpaper} that, for all $f \in \cD(A_{w,1})=\cD(A^{(w)})$,
\[
	\|B^{(w)}f\|_w \le \kappa \|A^{(w)}f\|_w \le \kappa \|A_{w,1}f\|_w.
\]
This, together with \eqref{C167}, implies
\[
	(1-\kappa)\|A_{w,1}f\|_w \le \|G_{w,1}f\|_w \le (1+\kappa)\|A_{w,1}f\|
\]
for all $f\in\cD(G_{w,1})=\cD(A_{w,1})$.
Since $\kappa<1$ by assumption, \eqref{C90} follows from \cite[Corollary~2.2.3]{lunardi1995analytic} 
and \eqref{C133}.
\end{proof}

\subsubsection{Proof of \cref{C100}}\label{C157}

We use \cref{C31,C32} in the proof of \cref{C100}.
Let us fix some notation.
In both cases, \ref{C128} and \ref{C129} in \cref{C98}, we take
\[
	X = \ell_w^1, \qquad Y = \ell_{\tildew}^1,
\]
which are Banach lattices, and $Y$ is a sublattice of $X$.
Note that $\tildew=w$ and hence $X=Y$ in case \ref{C128}.
Further, in case \ref{C128} we use
\[
	G = G^{(w)}, \qquad F(t,f) = K^{(\tildew,w)}(t,f), \qquad H = I;
\]
in case \ref{C129} we take
\[
	G = G_{w,1}, \qquad F(t,f) = K^{(\tildew,w)}(t,f) + f, \qquad H = (-A_{w,1})^\alpha.
\]
Note that also in case \ref{C128} we have $H=(-A_{w,1})^\alpha$ since $\alpha=0$ then.

Before we come to the proof of \cref{C100} we need one more lemma.

\begin{lemma}\label{C91}
Let $H$ be as above,
let $\tau\in[0,T)$, $r>0$ and $\gamma\ge c(\tau)r$.  Then
\begin{equation}\label{C92}
	K^{(\tildew,w)}(t,g)+\gamma Hg \ge0
\end{equation}
for all $t\in[0,\tau]$ and $g\in(\ell_{\tildew}^1)_+$ with $\|g\|_{\tildew}\le r$.
\end{lemma}

\begin{proof}
Let $t\in[0,\tau]$ and $g\in(\ell_{\tildew}^1)_+$ with $\|g\|_{\tildew}\le r$.
We use \eqref{C116} and the monotonicity of $(w_n)_{n=1}^\infty$ to obtain, for $n\in\NN$,
\begin{align*}
	& \bigl[K^{(\tildew,w)}(t,g)+\gamma Hg\bigr]_n
	= \bigl[K^{(\tildew,w)}(t,g)+\gamma(-A_{w,1})^\alpha g\bigr]_n
	\\[1ex]
	&= \frac{1}{2}\sum_{j=1}^{n-1}k_{n-j,j}(t)g_{n-j}g_j - g_n\sum_{j=1}^\infty k_{n,j}(t)g_j
	+ \gamma(a_n+1)^\alpha g_n
	\\[1ex]
	&\ge -g_nc(\tau)\sum_{j=1}^\infty\frac{\tildew_n\tildew_j}{w_{n+j}}g_j
	+ \gamma(a_n+1)^\alpha g_n
	= \biggl[-c(\tau)\sum_{j=1}^\infty \frac{w_n}{w_{n+j}}\tildew_j g_j + \gamma\biggr](a_n+1)^\alpha g_n
	\\[1ex]
	&\ge \bigl[-c(\tau)\|g\|_{\tildew}+\gamma\bigr](a_n+1)^\alpha g_n
	\ge \bigl[-c(\tau)r+\gamma\bigr](a_n+1)^\alpha g_n
	\ge 0,
\end{align*}
which proves \eqref{C92}.
\end{proof}

\begin{proof}[Proof of \cref{C100}]
We show that the conditions (a)--(f) of \cref{C31,C32} are satisfied.
Sometimes we have to distinguish the two cases \ref{C128} and \ref{C129}.

In case \ref{C128}, assumptions (a), (b), (d) and (e) are satisfied by 
Remarks~\ref{C61}\,(v) and \ref{C97}\,(iii).

In case \ref{C129} we obtain from the beginning of \cref{C156} and \cref{C89} that 
the operator $G_{w,1}$ generates an analytic semigroup with strictly negative growth bound
and that $Y=D_{G_{w,1}}(\alpha,1)=\ell_{\tildew}^1$.
Moreover, it is easy to see that $-H$ generates a substochastic semigroup.
Hence (a), (b), (d) and (e) hold by \cref{C86}.

In both cases it follows from \cref{C115,C146} that condition (c) in \cref{C31} is satisfied.
Moreover, condition (f) follows from \cref{C91} since the extra term in $F$
in \ref{C129} is positive.

Now the statements in (i)--(iii) follow immediately from \cref{C31,C32}.
Note that \eqref{C139} and \eqref{C140} are equivalent by \cref{C79} with $H=I$ and $\gamma=\pm1$.

Finally, we prove the claims in (iv).
Let $X_{[1]}$ and $\totmass$ be defined as in \eqref{C147} and \eqref{C159}, respectively,
and let $(S^{[1]}(t))_{t\ge0}$ be the semigroup in $X_{[1]}$, which exists
by \cite[Theorem~3.4]{kerr_lamb_langer2020fragpaper} since \eqref{C149} holds in both cases.
The space $\ell_w^1$ is continuously embedded as a dense subspace in $X_{[1]}$ by \cref{C88}, 
and $\totmass$, restricted to $\ell_w^1$, is a bounded linear functional on $\ell_w^1$.
Further, the semigroup $(S^{(w)}(t))_{t\ge0}$ extends to the semigroup $(S^{[1]}(t))_{t\ge0}$
since the two semigroups coincide on the set of sequences with only finitely many non-zero entries;
see also the discussion at the end of Section~3 in \cite{kerr_lamb_langer2020fragpaper}.
Moreover, \cref{C134} with $\omega_n=n$ implies that
\begin{equation}\label{C190}
	\totmass\bigl(K^{(\tildew,w)}(t,f)\bigr) = 0
\end{equation}
for all $t\in[0,T)$ and $f\in\ell_{\tildew}^1$.

Let us first consider the case when \eqref{C149} holds.
It follows from \cite[Theorem~3.4]{kerr_lamb_langer2020fragpaper} that 
the semigroup $(S^{[1]}(t))_{t\ge0}$ is substochastic, 
and hence \eqref{C176} is satisfied with $\phi=\totmass$.
Relation \eqref{C190} shows that \eqref{C177} holds with $\nu(\tau)=0$.
Hence we can apply \cref{C143}, which yields \eqref{C165}.

Assume now that \eqref{C26} holds.
It follows again from \cite[Theorem~3.4]{kerr_lamb_langer2020fragpaper} that 
the semigroup $(S^{[1]}(t))_{t\ge0}$ is stochastic.
Together with \cite[Proposition~2.3\,(i)]{kerr_lamb_langer2020fragpaper}, 
this shows that \eqref{C123} is satisfied.
Clearly, \eqref{C124} follows from \eqref{C190}.
Hence, we can apply \cref{C84}, which proves \eqref{C189}.
\end{proof}

\subsection{Classical Solutions}\label{C158}

In this section we show that, under an additional differentiability assumption on the coagulation rates, 
the mild solution in \cref{C100} is also a classical solution.
We deal with the two cases, (CI) and (CII), separately.  Let us start with the former.

\begin{theorem}\label{C101}
Let \cref{C22,C88,C98} with \ref{C128} hold and assume that $k_{n,j}$ are continuously
differentiable for all $n,j\in\NN$ and that, for every $t'\in(0,T)$, there 
exists $\tilde c(t')>0$ such that
\begin{equation}\label{C163}
	|k'_{n,j}(t)| \le \tilde c(t')\frac{w_n w_j}{w_{n+j}}
	\qquad \text{for all} \ n,j\in\NN,\;t\in[0,t'].
\end{equation}
Further, let $\mr u\in(\cD(G^{(w)}))_+$ and let $u$ be the unique mild solution
of \eqref{C137}, \eqref{C138}.  Then $u$ is also a classical solution of \eqref{C137}, \eqref{C138}.
\end{theorem}

\begin{proof}
Let $\widetilde K^{(w,w)}$ be as in \eqref{C162}, define
\[
	\widetilde L^{(w,w)}[t,f,g] \coloneq \biggl(\frac{1}{2}\sum_{j=1}^{n-1} k'_{n-j,j}(t)f_{n-j}g_j
	- \sum_{j=1}^{\infty} k'_{n,j}(t)f_n g_j\biggr)_{n=1}^\infty,
	\quad f,g\in\ell_w^1,\; t\in(0,T),
\]
and set
\[
	h_{n,j}(t) \coloneq \frac{w_{n+j}}{w_n w_j}k_{n,j}(t), \qquad n,j\in\NN,\; t\in[0,T).
\]
We show that the assumptions of \cref{C46} are satisfied.
Let $t_0\in(0,T)$ and choose $t'\in(t_0,T)$.  For $f,g\in\ell_w^1$ and $t\in[0,t']$ we obtain from \cref{C110}
with 
\[
	\varphi_{n,j}=\frac{k_{n,j}(t)-k_{n,j}(t_0)}{t-t_0}-k'_{n,j}(t_0)
\]
that
\begin{align}
	& \bigg\|\frac{1}{t-t_0}\Bigl(\widetilde K^{(w,w)}[t,f,g]-\widetilde K^{(w,w)}[t_0,f,g]\Bigr)
	-\widetilde L^{(w,w)}[t_0,f,g]\bigg\|_w
	\nonumber\\[1ex]
	&\le \frac{3}{2}\sum_{n=1}^\infty\sum_{j=1}^\infty 
	w_{n+j}\bigg|\frac{k_{n,j}(t)-k_{n,j}(t_0)}{t-t_0}-k'_{n,j}(t_0)\bigg|\,|f_n|\,|g_j|
	\nonumber\\[1ex]
	&= \frac{3}{2}\sum_{n=1}^\infty\sum_{j=1}^\infty \bigg|\frac{h_{n,j}(t)-h_{n,j}(t_0)}{t-t_0}-h'_{n,j}(t_0)\bigg|
	\cdot w_n w_j|f_n|\,|g_j|.
	\label{C164}
\end{align}
It follows from \eqref{C163} and the Mean Value Theorem that the expression in the series in \eqref{C164}
is bounded by $2\tilde c(t')w_nw_j|f_n|\,|g_j|$.  The differentiability of $k_{n,j}$ implies that
the expressions in the series in \eqref{C164} converge to $0$ as $t\to t_0$ for all $n,j\in\NN$.
Hence the double series in \eqref{C164} converges to $0$ as $t\to t_0$ by the Dominated Convergence Theorem.
This shows that $\frac{\partial}{\partial t}\widetilde K^{(w,w)}[t,f,g]=\widetilde L^{(w,w)}[t,f,g]$
for all $f,g\in\ell_w^1$.  We can again use \cref{C110}, now with $\varphi_{n,j}=k'_{n,j}(t)$,
to obtain
\[
	\big\|\widetilde L^{(w,w)}[t,f,g]\big\|_w 
	\le \frac{3}{2}\sum_{n=1}^\infty\sum_{j=1}^\infty \big|h'_{n,j}(t)\big|\cdot w_n w_j|f_n|\,|g_j|
	\le \frac{3}{2}\tilde c(t')\|f\|_w\|g\|_w
\]
for $t\in(0,t']$ and $f,g\in\ell_w^1$.
In a similar way as above one shows the continuity of $\frac{\partial}{\partial t}\widetilde K^{(w,w)}[t,f,g]$
by applying again \cref{C110}, the continuity of $k'_{n,j}$ and the Dominated Convergence Theorem.
Hence all assumptions of \cref{C46} are satisfied; the latter implies that $K^{(w,w)}$ is 
Fr\'echet differentiable with continuous Fr\'echet derivative.
Now the claim follows from \cite[Theorem~6.1.5]{pazy1983semigroups}.
\end{proof}

\noindent
The following theorem deals with classical solutions in the case (CII).

\begin{theorem}\label{C102}
Let \cref{C22,C88,C98} with \ref{C129} hold and assume that, for every $t'\in(0,T)$, 
there exist $\tilde c(t')>0$ and $\theta(t')\in(0,1)$ such that
\begin{equation}\label{C166}
	|k_{n,j}(t)-k_{n,j}(s)| \le \tilde c(t')\frac{\tildew_n\tildew_j}{w_{n+j}}(t-s)^{\theta(t')}
	\qquad\text{for all} \ n,j\in\NN,\; s,t\in[0,t'].
\end{equation}
Further, let $\mr u\in(\ell_{\tildew}^1)_+$ and let $u$ be the unique mild solution
of \eqref{C137}, \eqref{C138}.  Then $u$ is also a classical solution of \eqref{C137}, \eqref{C138}.
\end{theorem}

\begin{proof}
Let $t'\in(0,T)$.
We obtain from \cref{C110} with $\varphi_{n,j}=k_{n,j}(t)-k_{n,j}(s)$ and \eqref{C166} that, 
for $f\in\ell_{\tildew}^1$ and $s,t\in[0,t']$,
\begin{align*}
	\big\|K^{(\tildew,w)}(t,f)-K^{(\tildew,w)}(s,f)\big\|_w
	&\le \frac{3}{2}\sum_{n=1}^\infty\sum_{j=1}^\infty w_{n+j}|k_{n,j}(t)-k_{n,j}(s)|\cdot|f_n|\,|f_j|
	\\[1ex]
	&\le \frac{3}{2}\tilde c(t')\sum_{n=1}^\infty\sum_{j=1}^\infty (t-s)^{\theta(t')}\tildew_n \tildew_j|f_n|\,|f_j|
	\\[1ex]
	&= \frac{3}{2}\tilde c(t')(t-s)^{\theta(t')}\|f\|_{\tildew}^2.
\end{align*}
It follows from \cite[Theorem~7.1.10]{lunardi1995analytic} that $u$ is a classical solution on $[0,t']$.
Since $t'$ was arbitrary in $(0,T)$, the theorem is proved.
\end{proof}

\subsection{A Global Existence Result}\label{C152}

Finally, in the setting (CI), we show that, under an extra assumption on the weight and the coagulation coefficients,
we have $\tmax=T$, i.e.\ we show that the unique mild solution in Theorem~\ref{C100} is a global solution.

\begin{theorem}\label{C141}
Let \cref{C22,C88,C98} with \ref{C128} be satisfied where $w_n=n^p$ with some $p>1$ and assume 
that \eqref{C149} holds.
Let $\mu:(0,T)\to(0,\infty)$ be monotonic increasing function such that,
for each $t'\in(0,T)$,
\begin{equation}\label{C191}
	k_{n,j}(t) \le \mu(t')(n+j) \qquad\text{for all} \ n,j\in\NN, \, t\in[0,t'].
\end{equation}
Let $\mr u\in(\ell_w^1)_+$ and let $u$ be the unique mild solution of \eqref{C137}, \eqref{C138} on $[0,\tmax)$.
Then $\tmax=T$, i.e.\ $u$ is a global solution.
Moreover,
\begin{equation}\label{C144}
	\|u(t)\|_w \le \|\mr u\|_w e^{c_p\mu(t)\totmass(\mr u)t}, \qquad t\in[0,T),
\end{equation}
with
\begin{equation}\label{C145}
	c_p = \begin{cases}
		p & \text{when} \ p\in(1,2],
		\\[1ex]
		2^p-2 & \text{when} \ p\in(2,\infty).
	\end{cases}
\end{equation}
\end{theorem}

\begin{proof}
First note that \cref{C100} implies that $u(t)\ge0$ for all $t\in[0,\tmax)$ and that \eqref{C165} holds.
It follows from \cite[Lemma~2.3]{carr1992asymptotic} 
(see also \cite[Lemma~7.4.4]{banasiak_lamb_laurencot2020vol2}) that
\begin{equation}\label{C142}
	(x+y)\bigl[(x+y)^p-x^p-y^p\bigr] \le c_p(x^py+xy^p), \qquad x,y\in(0,\infty),
\end{equation}
where $c_p$ is as in \eqref{C145}.
Let $\phi_w$ be the unique linear extension of $\normcdotsub{\ell_w^1}$ from the positive cone
to $\ell_w^1$, i.e.\ $\phi_w(f) = \sum_{n=1}^\infty w_nf_n$.
Now let $t'\in(0,\tmax)$.
Using \eqref{C142}, \cref{C134} with $\omega_n=w_n$, and \eqref{C165} we obtain, for $t\in[0,t']$,
\begin{align}
	\phi_w\bigl(K^{(\tildew,w)}(t,u(t))\bigr)
	&= \frac{1}{2}\sum_{n=1}^\infty\sum_{j=1}^\infty \bigl(w_{n+j}-w_n-w_j\bigr)k_{n,j}(t)u_n(t)u_j(t)
	\nonumber\\[1ex]
	&\le \frac{1}{2}\sum_{n=1}^\infty\sum_{j=1}^\infty \bigl[(n+j)^p-n^p-j^p\bigr]\mu(t')(n+j)u_n(t)u_j(t)
	\nonumber\\[1ex]
	&\le \frac{1}{2}c_p\mu(t')\sum_{n=1}^\infty\sum_{j=1}^\infty \bigl(n^p j+nj^p\bigr)u_n(t)u_j(t)
	\nonumber\\[1ex]
	&= c_p\mu(t')\sum_{j=1}^\infty ju_j(t)\cdot \sum_{n=1}^\infty n^pu_n(t)
	\nonumber\\[1ex]
	&\le c_p\mu(t')\totmass(\mr u)\phi_w\bigl(u(t)\bigr).
	\label{C160}
\end{align}
Since $(S^{(w)}(t))_{t\ge0}$ is substochastic on $\ell_w^1$,
we can again apply \cref{C143} with $\phi=\phi_w$, which yields
\begin{equation}\label{C93}
	\phi_w\bigl(u(t)\bigr) \le \phi_w(\mr u)e^{c_p\mu(t')\totmass(\mr u)t}, \qquad t\in[0,t'].
\end{equation}
If $\tmax<T$, then $\|u(t)\|_w\le\|\mr u\|_w e^{c_p\mu(\tmax)\totmass(\mr u)\tmax}$ for all $t\in[0,\tmax)$,
which contradicts \cref{C100}\,(ii).
Hence $\tmax=T$, and \eqref{C144} follows from \eqref{C93}.
\end{proof}

\begin{remark}
Note that, to obtain the global existence result in Theorem~\ref{C141}, we require a stronger assumption 
on the coagulation rates than was earlier imposed in \eqref{C116}.  For example, if we consider the diagonal $n=j$, 
then, setting $\alpha=0$ and $w_n=n^p$ for some $p>1$, \eqref{C116} becomes $k_{n,n}(t) \le \frac{c(t')}{2}n^p$, 
while \eqref{C191} 
is $k_{n,n}(t) \leq 2\mu(t')n$. 
Hence, whereas \eqref{C116} allows power-law growth in the coagulation rates along the diagonal, 
global existence requires the diagonal to grow no more than linearly.
\end{remark}

\noindent
\textbf{Acknowledgements.} \\
L.~Kerr gratefully acknowledges the support of \emph{The Carnegie Trust for the Universities of Scotland}.
Moreover, both authors thank Wilson Lamb for many fruitful discussions.
\\[1ex]
\textbf{Data Availability:} 
The manuscript has no associated data.
\\[1ex]
\textbf{Financial or non-financial interest:}
The authors have no relevant financial or non-financial interests to disclose.


%

\end{document}